\RequirePackage{fix-cm}
\let\accentvec\vec
\documentclass{svjour3}

\smartqed

\let\vec\accentvec

\usepackage{graphicx}
\graphicspath{{./figures/}}
\usepackage{amssymb,amsfonts,amsmath}
\usepackage{bm}
\usepackage[caption=false]{subfig}
\usepackage{color}
\usepackage{dashrule}
\usepackage{setspace}
\usepackage{url}
\usepackage[numbers,sort&compress]{natbib}

\newcommand{\beq}{\begin{equation}}
\newcommand{\eeq}{\end{equation}}
\newcommand{\lbracket}{[\![}
\newcommand{\rbracket}{]\!]}

\journalname{}

\begin{document}

\title{A Simple and Efficient Preconditioning Scheme for Heaviside Enriched XFEM}


\author{Christapher Lang \and David Makhija \and Alireza Doostan \and Kurt Maute}

\institute{C. Lang \at Structural Mechanics and Concepts Branch, NASA Langley Research Center, Hampton, VA \and D. Makhija \and A. Doostan \and K. Maute \at Aerospace Engineering Sciences, University of Colorado, Boulder, CO \\
\email{alireza.doostan@colorado.edu}}

\date{Received: date / Accepted: date}

\maketitle

\begin{abstract}
The eXtended Finite Element Method (XFEM) is an approach for solving problems with non-smooth solutions, which arise from geometric features such as cracks, holes, and material inclusions. In the XFEM, the approximate solution is locally enriched to capture the discontinuities without requiring a mesh which conforms to the geometric features. One drawback of the XFEM is that an ill-conditioned system of equations results when the ratio of volumes on either side of the interface in an element is small. Such interface configurations are often unavoidable, in particular for moving interface problems on fixed meshes. In general, the ill-conditioning reduces the performance of iterative linear solvers and impedes the convergence of solvers for nonlinear problems. This paper studies the XFEM with a Heaviside enrichment strategy for solving problems with stationary and moving material interfaces. A generalized formulation of the XFEM is combined with the level set method to implicitly define the embedded interface geometry. In order to avoid the ill-conditioning, a simple and efficient scheme based on a geometric preconditioner and constraining degrees of freedom to zero for small intersections is proposed. The geometric preconditioner is computed from the nodal basis functions, and therefore may be constructed prior to building the system of equations. This feature and the low-cost of constructing the preconditioning matrix makes it well suited for nonlinear problems with fixed and moving interfaces. It is shown by numerical examples that the proposed preconditioning scheme performs well for $C^0$-continuous problems with both the stabilized Lagrange and Nitsche methods for enforcing the continuity constraint at the interface. Numerical examples are presented which compare the condition number and solution error with and without the proposed preconditioning scheme. The results suggest that the proposed preconditioning scheme leads to condition numbers similar to that of a body-fitted mesh using the traditional finite element method without loss of solution accuracy.


\end{abstract}

\section{Introduction}

A standard tool for numerically solving problems defined by a set of partial differential equations in many engineering disciplines is the Finite Element Method (FEM). The solution to problems which feature embedded interfaces, such as material inclusions or voids, is non-smooth due to strong or weak discontinuities which occur at the interface. A strong discontinuity occurs when the solution is discontinuous across the interface. A weak discontinuity occurs when the solution is continuous but its spatial derivatives are discontinuous across the interface. Conventionally, a finite element mesh is used which conforms to the interface in order to approximate the non-smooth solution. However, mesh generation may lead to robustness issues and increase the computational cost for problems with complex geometries or moving interfaces.

\begin{figure}[ht]
\centering
\subfloat[]{\includegraphics[trim=1.2in 5in 1.2in 1.8in,clip,width=0.45\textwidth]{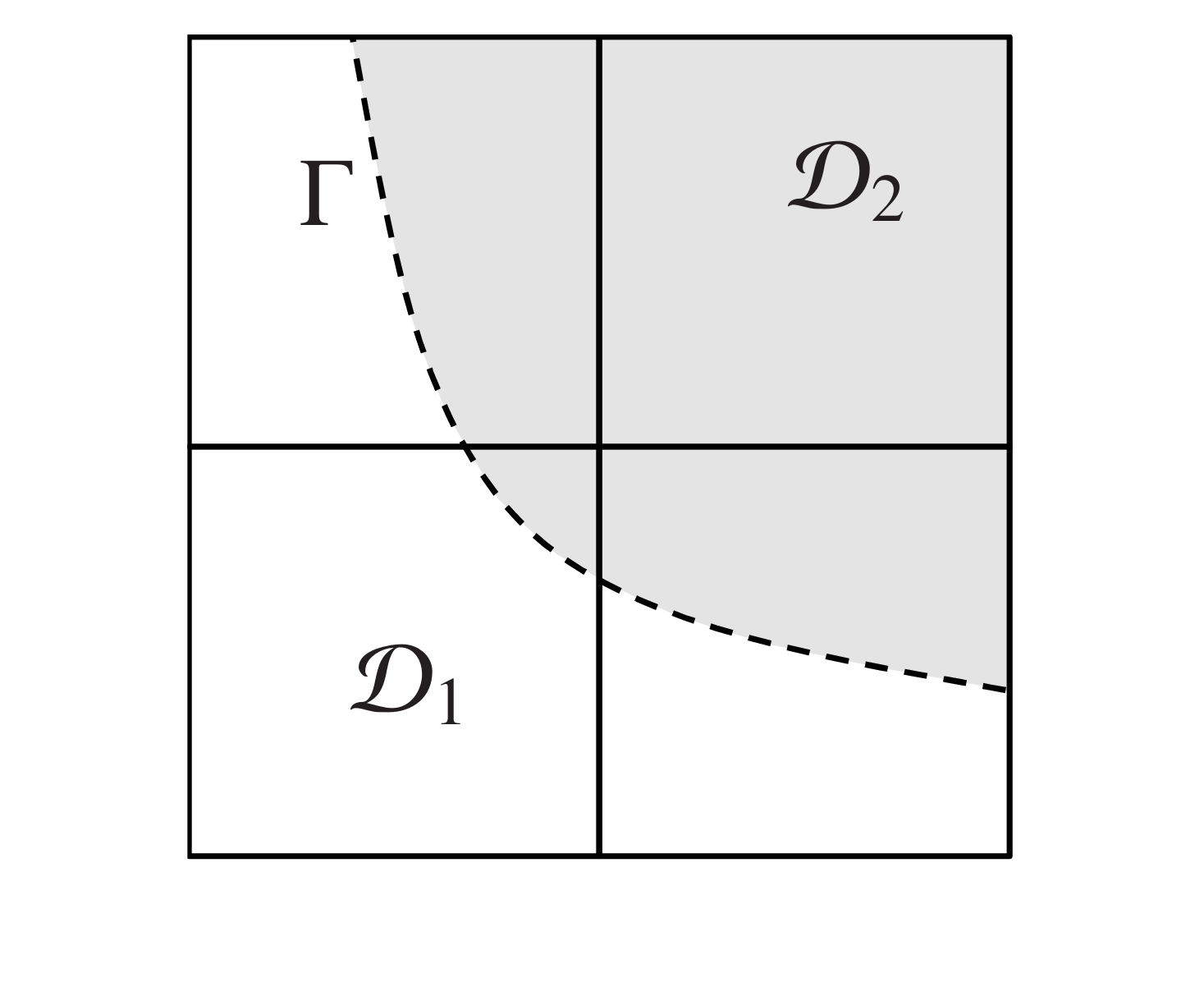}}
\subfloat[]{\includegraphics[trim=1.2in 5in 1.2in 1.8in,clip,width=0.45\textwidth]{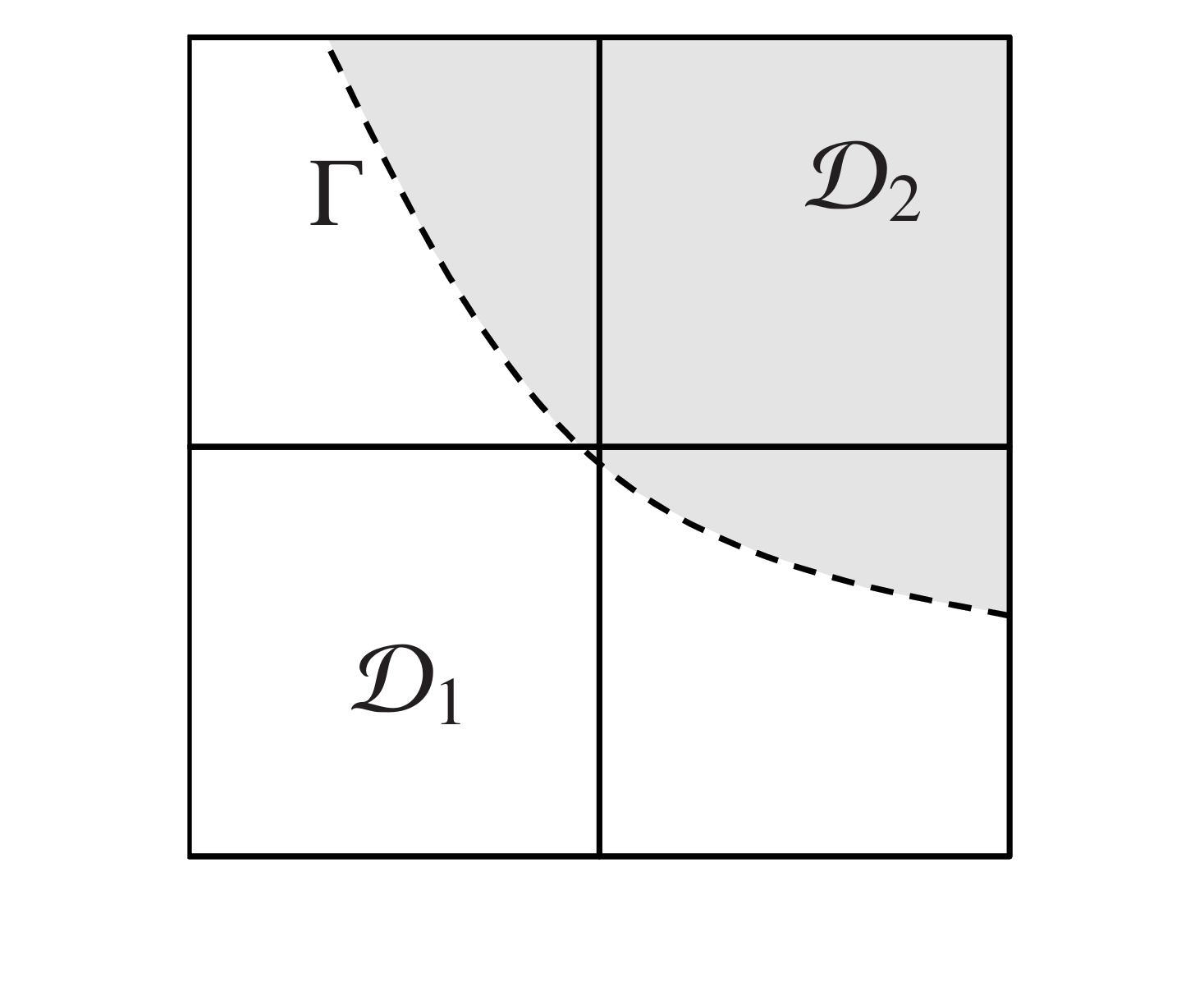}}
\caption{Configuration of four quadrilateral elements which lead to a (a) well-conditioned and (b) ill-conditioned system. The lower left element in (b) has a small ratio of areas bisected by the interface.}
\label{fig:configs}
\end{figure}

A widely used alternative for solving problems with embedded interfaces is the eXtended Finite Element Method (XFEM) \citep{ref:moes1999,ref:sukumar2001}. Local enrichment functions are added to the standard FEM basis to represent the strong or weak solution discontinuities. The enrichment functions are constructed based on the position of the interface, which is implicitly defined by the level set method \citep{ref:osher1988,ref:sethian1999}. The XFEM does not require a mesh that conforms to the interface, which reduces the complexity of mesh construction. This feature is particularly advantageous for complex geometries as well as problems with moving or changing interface configurations \citep{ref:chessa2002,ref:zabaras2006,KM:12,ref:lang2012}. However, the XFEM can lead to ill-conditioned systems when an intersected element(s) has a small ratio of areas bisected by the interface, as illustrated in Fig. \ref{fig:configs}. Ill-conditioned systems are a particular issue for nonlinear problems and iterative linear solvers \citep{ref:bechet2005,ref:fries2010}.

The focus of this work is on a new scheme to mitigate the ill-conditioning issue in the XFEM. The goal is to obtain condition numbers using the XFEM that are of the same order of magnitude as standard FEM with a conforming mesh. Various approaches for dealing with this ill-conditioning have been proposed. A straight-forward approach is to construct a mesh that avoids small intersections with a uniform ratio of intersected element areas. Another approach is to move the nodes of intersected elements in order to avoid any intersected areas less than a specified amount \citep{ref:choi2012}. However these approaches require adaptive meshing and mesh updating strategies which typically encounter efficiency and robustness issues for complex geometries and moving interfaces. 

Other approaches involve modifications to the discretized system of equations such that careful mesh construction or moving the nodes is not necessary. Reusken \citep{ref:reusken2008} suggested constraining degrees of freedom associated with small supports to zero. This approach improves the condition number of the system by removing the constrained degrees of freedom. However, there is a trade-off between the accuracy of the solution and the ill-conditioning of the system which depends on the criteria for selecting the degrees of freedom to be constrained. The criteria must be carefully chosen in order to improve the condition number without decreasing the solution accuracy beyond an acceptable level. Preconditioning schemes have been proposed to improve the condition number of the system matrices to be solved. Sauerland and Fries \citep{ref:sauerland2012} study a Jacobi preconditioner, and preconditioners based on a Cholesky decomposition are studied by Bechet et al \citep{ref:bechet2005} and Menk and Bordas \citep{ref:menk2011}. These alternative schemes are well suited for linear problems. However, the preconditioner can be built only after the discretized system of equations is assembled and must be reconstructed in each solution step for nonlinear problems, even when the interface geometry remains fixed.

A third class of methods modify the enrichment function to avoid the ill-conditioning issue. In \citep{ref:ruberg2012}, an approach for dealing with small intersections using b-spline finite elements is introduced. Interior and exterior b-splines are defined by the intersection size, and b-splines with a small intersection are denoted as exterior. The degrees of freedom associated with the exterior b-splines are expressed by a linear combination of the interior b-splines degrees of freedom. A stable XFEM is described in \citep{ref:babuska2012,ref:sauerland2012} which uses a local enrichment function constructed from a linear interpolant of the global enrichment function in the intersected elements.

Finally, Hansbo et al \citep{ref:hansbo2013} and Wadbro et al \citep{ref:wadbro2013} propose to augment the weak formulation to produce a well-conditioned system of equations independent of the interface position. The solution for each subdomain separated by the interface is considered, and a version of Nitsche's method is used to enforce the interface conditions. By adding additional volume terms to the weak formulation, the ill-conditioning is mitigated, but the solution error at the interface is increased. While this error decreases with mesh refinement, for a given mesh size this approach alters the solution of the discretized system.

In this work, a preconditioning scheme is proposed for a generalized Heaviside enrichment \citep{ref:makhija2013} that consists of a linear preconditioner and constraining degrees of freedom associated with small intersections. For the proposed scheme, no special considerations are necessary in the mesh generation, the enrichment function is not modified, and the weak formulation is unchanged. The construction of the preconditioner only requires the nodal basis functions and interface geometry; therefore, it may be constructed prior to building the discretized system of equations and is well suited for nonlinear problems.

Problems with static and prescribed moving interfaces are studied, and numerical examples show condition numbers for the XFEM using the proposed preconditioning scheme similar to the standard FEM. The proposed approach shows satisfactory performance for the stabilized Lagrange and Nitsche methods \citep{ref:stenberg1995,ref:mendez2004,ref:juntunen2008} for enforcing continuity at the interface.

The remainder of this paper is organized as follows: Section \ref{sec:model} defines the model problem for this work. Section \ref{sec:xfem} describes the XFEM framework, Heaviside enrichment strategy, and interface constraint formulation. Section \ref{sec:precon} presents the proposed preconditioning scheme for handling small intersections. In Section \ref{sec:examples}, three numerical examples are presented to demonstrate the key features of the projection scheme.

\section{Model Setup} \label{sec:model}

\begin{figure}[ht]
\centering\includegraphics[trim=1.8in 4in 0.9in 1.6in,clip,width=0.45\textwidth]{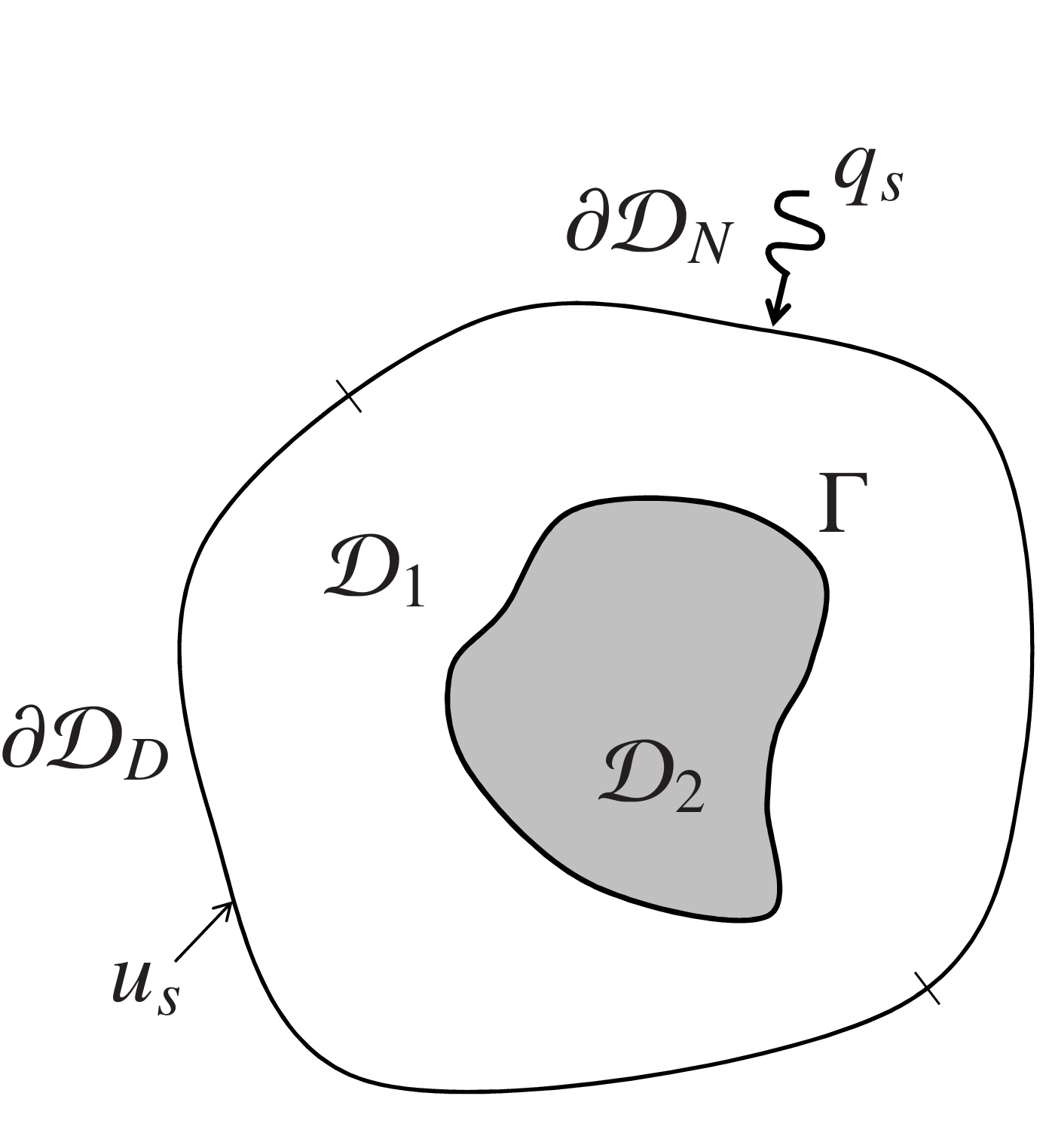}
\caption{Schematic of the model diffusion problem.}
\label{fig:model}
\end{figure}

Here we consider solving a stationary diffusion equation for a material with a single inclusion, as depicted in Fig. \ref{fig:model}. The model problem is used for the description of the numerical method and for the first two numerical examples of Section \ref{sec:examples}. While we focus on this model problem for describing the details of the preconditioning scheme, the method is applicable to other problem types. In particular, the performance of the preconditioning scheme for a transient nonlinear fluid flow problem with moving interfaces is presented in the third numerical example in Section \ref{sec:examples}.

The domain is comprised of two non-overlapping subdomains, such that $\mathcal{D}=\mathcal{D}_1 \cup \mathcal{D}_2$ and $\mathcal{D}_1\cap\mathcal{D}_2=\emptyset$. The interface between the two subdomains is defined as $\Gamma = \partial \mathcal{D}_1 \cap \partial \mathcal{D}_2$. A level set function $\phi(\bm{x})$ is constructed to define the location of $\Gamma$, such that 
\begin{align}
& \phi(\bm{x}) < 0 \quad \text{if $\bm{x} \in \mathcal{D}_1$} \nonumber \\
& \phi(\bm{x}) > 0 \quad \text{if $\bm{x} \in \mathcal{D}_2$} \nonumber \\
& \phi(\bm{x}) = 0 \quad \text{if $\bm{x} \in \Gamma$ .}
\end{align}

\noindent In this work, the signed distance function is used to define the level set function,
\beq
\phi(\bm{x}) = \pm \min \| \bm{x} - \bm{x}_\Gamma \| \text{ ,}
\label{eq:signdist}
\eeq
\noindent where $\bm{x}_\Gamma$ is the interface location and $\| \cdot \|$ denotes the $L^2$-distance. Considering the particular case of diffusive heat conduction, the model problem consists of finding the temperature distribution, $u(\bm{x})$, such that 
\begin{align}
- \nabla \cdot (\bm{\kappa} \nabla u_i) & = f \quad \text{ in $\mathcal{D}_i$} \nonumber \\
u_i & =u_s \quad \text{on $\partial\mathcal{D}_i \cap \partial\mathcal{D}_D$} \nonumber \\
(\bm{\kappa} \nabla u_i) \cdot \bm{n}_i & =q_s \quad \text{on $\partial\mathcal{D}_i \cap \partial\mathcal{D}_N$}
\label{eq:model}
\end{align}
\noindent for $i=1,2$, where $\bm{\kappa}$ is the thermal conductivity tensor, $f$ is a volumetric heat source, and $u_i$ denotes the restriction of $u$ to $\mathcal{D}_i$. The temperature distribution $u_s$ is specified on a Dirichlet boundary $\partial\mathcal{D}_D$, and the heat flux $q_s$ is specified on a Neumann boundary $\partial\mathcal{D}_N$. The outward unit normal to $\mathcal{D}_i$ is denoted by $\bm{n}_i$. Additionally, continuity of the solution and flux across the interface $\Gamma$ must be satisfied, such that
\begin{align}
& \lbracket u \rbracket = u_1-u_2 = 0 \quad \text{on $\Gamma$} \nonumber\\
& k_1\nabla u_1 \cdot \bm{n}_1 + k_2\nabla u_2 \cdot \bm{n}_2 = 0 \quad \text{on $\Gamma$ .}
\label{eq:ucont}
\end{align}

Without loss of generality, the materials are assumed to be isotropic, i.e. $\bm{\kappa}=k \ \mathbf{I}$. The conductivity $k$ is defined as
\beq
k(\bm{x})=\left\{
\begin{array}{l l}
k_1 \quad \text{if $\bm{x} \in \mathcal{D}_1$} \\
k_2 \quad \text{if $\bm{x} \in \mathcal{D}_2$}
\end{array}
\right.
\eeq
\noindent with constants $k_1$ and $k_2$.

\section{Extended Finite Element Method} \label{sec:xfem}

The traditional finite element method requires a mesh which conforms to the interface to implicitly satisfy the temperature continuity and to capture the discontinuity in the temperature gradients at $\Gamma$. Alternatively, the extended finite element method is used to locally capture the non-smooth solution at the interface without using a mesh which conforms to $\Gamma$. In this section, we briefly outline the particular XFEM approach used here for solving the governing equation in (\ref{eq:model}).

The weak form of the governing equations is constructed by multiplying (\ref{eq:model}) by a set of admissible test functions and integrating over $\mathcal{D}$. The space $V=H^1(\mathcal{D})$ is the Hilbert space consisting of functions with square integrable first derivatives and $V_0 = \{v \in V:v|_{\partial\mathcal{D}_D} = 0\}$. Let $u \in V$ be the solution and $v \in V_0$ be an admissible test function.  The weak form of the model problem is stated as: Find $u \in V$ such that $u=u_s$ on $\partial\mathcal{D}_D$ and
\beq
\int_{\mathcal{D}} (\bm{\kappa} \nabla u) \cdot \nabla v \ d\bm{x} - \int_{\mathcal{D}} fv \ d\bm{x} - \int_{\partial\mathcal{D}_N} q_s v \ ds =0 \quad \forall v \in V_0 \text{ .}
\label{eq:weak}
\eeq
\noindent Note that the continuity conditions were used to express the weak form in (\ref{eq:weak}), such that $\lbracket u \rbracket = 0$ at $\Gamma$ and
\beq
\int_{\Gamma} (k_1\nabla u_1 \cdot \bm{n}_1) v \ ds + \int_{\Gamma} (k_2\nabla u_2 \cdot \bm{n}_2) v \ ds =0 \text{ .}
\label{eq:bterms}
\eeq

In the XFEM, the traditional finite element approximation is augmented by an enrichment function and additional degrees of freedom. The choice of enrichment function affects the convergence and accuracy of the approximation, and various types of enrichment functions have been proposed. A $C^0$-continuous enrichment function \citep{ref:moes2003} inherently satisfies the solution continuity at $\Gamma$. As discussed in \citep{ref:fries2008}, the nodes of neighboring elements to intersected elements, called blending elements, also require enriched degrees of freedom for accurate solutions. A step enrichment function, such as a Heaviside or sign function, simplifies the formulation since enriched nodal basis functions and blending elements are not required. However, the approximation of the weak form (\ref{eq:weak}) needs to be augmented by constraints to satisfy the temperature continuity at the interface. Both $C^0$-continuous and step enrichment functions can lead to a system of equations that is ill-conditioned \citep{ref:soghrati2010,ref:soghrati2012}.

Here, we follow the work of Terada et al \citep{ref:terada2003} and adopt a generalized version of the Heaviside enrichment strategy of Hansbo and Hansbo \citep{ref:hansbo2004}. As recently shown by Makhija and Maute \citep{ref:makhija2013}, this implementation of the XFEM provides great flexibility in discretizing a broad range of partial differential equations with multiple phases for any choice of nodal basis functions. The remainder of this section describes the details of the generalized Heaviside enrichment strategy and the interface constraint formulation.

\subsection{Generalized Heaviside Enrichment} \label{sec:heav}

Consider a finite element mesh, $\mathcal{T}_h$, for $\mathcal{D}$ consisting of elements with edges that do not necessarily coincide with $\Gamma$. A Heaviside enrichment function is implemented in the XFEM formulation such that the approximation to the solution for two phases is defined as
\begin{align}
\hat{u}(\bm{x}) =\sum_{m=1}^{M} &\left(H(-\phi(\bm{x}))\sum_{i \in I} N_i(\bm{x}) u_{i,m}^{(1)} \right. \nonumber\\
& \quad \left. + H(\phi(\bm{x}))\sum_{i \in I} N_i(\bm{x}) u_{i,m}^{(2)}  \right)
\label{eq:genxfem}
\end{align}
\noindent where $I$ is the set of all nodes in $\mathcal{T}_h$, $N_i(\bm{x})$ are the nodal basis functions, $M$ is the maximum number of enrichment levels, $u_{i,m}^{(p)}$ is the degree of freedom at node $i$ for phase $p\in\{1,2\}$, and $H$ is the Heaviside function,
\beq
H(z)=\left\{
  \begin{array}{l l}
  1 & \quad z>0 \\
  0 & \quad z\le 0
  \end{array}
  \right. \text{ .}
  \label{eq:heaviside}
\eeq

\begin{figure}[ht]
\centering\includegraphics[trim=1.2in 3.5in 0.0in 1.9in,clip,width=0.45\textwidth]{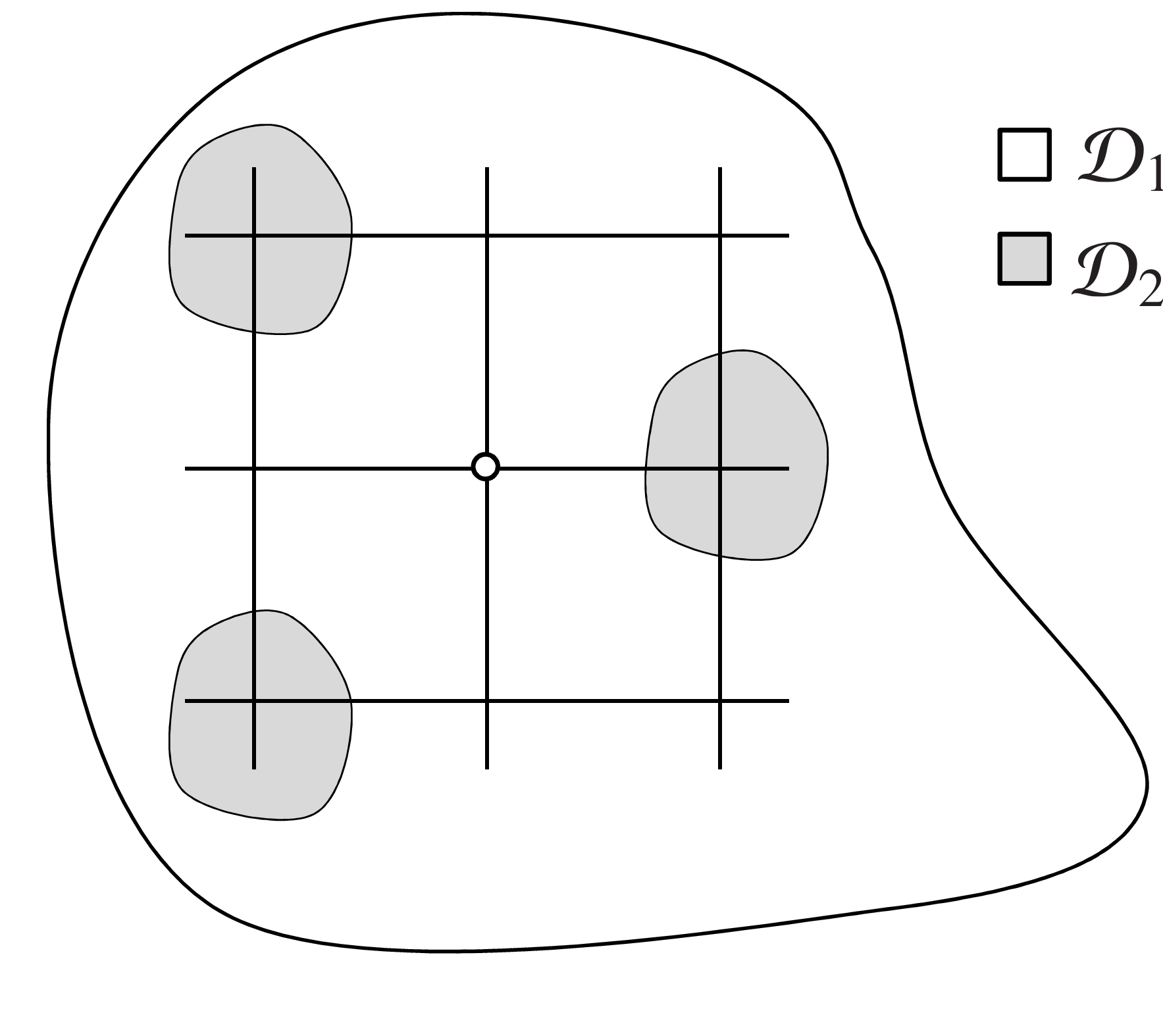}
\caption{Example configuration requiring multiple enrichment levels at the center node.}
\label{fig:config3}
\end{figure}

The need for multiple enrichment levels is illustrated by the example configuration shown in Fig. \ref{fig:config3}. Four quadrilateral elements share a central node that is connected to the phase 1 domain and three inclusions belonging to phase 2. The center node requires one degree of freedom for the phase 1 solution and three degrees of freedom in order to individually interpolate the solutions in the three inclusions. By generalizing the Heaviside enrichment to multiple levels, accurate solutions can be determined for neighboring intersected elements and elements intersected more than once. The number of enrichment levels required at a single node is determined by the number of disconnected regions of the same phase included in the support of the nodal basis function. Note that while a maximum number of enrichment levels is specified in (\ref{eq:genxfem}), some enrichment levels are not used. The degrees of freedom corresponding to the unused enrichment levels are removed from the system of equations. Further details of this generalized enrichment strategy is provided in \citep{ref:makhija2013}.

\vspace{-1cm}
\subsection{Interface Constraint Formulation} \label{sec:ifc}

While the continuity of the solution at the interface \eqref{eq:ucont} is inherently satisfied using a $C^0$-continuous enrichment function, the Heaviside enrichment requires an additional constraint to enforce the continuity. Common approaches for enforcing an interface constraint in the XFEM include the stabilized Lagrange multiplier and Nitsche methods \citep{ref:stenberg1995,ref:mendez2004,ref:juntunen2008}. Here, both constraint formulations are used for enforcing continuity at the interface for the model problem.

The weak form using the stabilized Lagrange multiplier method is stated as: Find $(u,\lambda) \in (V \times W)$ such that $u=u_s$ on $\partial\mathcal{D}_D$ and 
\begin{align}
&\int_{\mathcal{D}} (\bm{\kappa} \nabla u) \cdot \nabla v d\bm{x} - \int_{\mathcal{D}} fv d\bm{x} - \int_{\partial\mathcal{D}_N} q_s v ds \nonumber\\
& \quad - \int_\Gamma \lbracket v \rbracket \lambda d\Gamma + \int_\Gamma \mu \left ( \lambda - \left \{ k \nabla u \cdot \bm{n} \right \} \right ) d\Gamma \nonumber\\
& \quad - \gamma_S \int_\Gamma \mu \lbracket u \rbracket d\Gamma =0 \quad \forall (v,\mu) \in (V_0 \times W) \text{ ,}
\label{eq:slag}
\end{align}
\noindent where $\lambda$ is the Lagrange multiplier, $W=H^{-1/2}(\Gamma)$ is the space for the Lagrange multiplier, $\mu$ is the associated test function, $\gamma_S$ is a constraint factor, and $\{\cdot\}=\frac{1}{2}(\cdot)_1+\frac{1}{2}(\cdot)_2$ denotes the mean operator on the interface.

For Nitsche's method, the weak form is stated as: Find $u \in V$ such that $u=u_s$ on $\partial\mathcal{D}_D$ and
\begin{align}
&\int_{\mathcal{D}} (\bm{\kappa} \nabla u) \cdot \nabla v d\bm{x} - \int_{\mathcal{D}} fv d\bm{x} - \int_{\partial\mathcal{D}_N} q_s v ds \nonumber\\ & \quad - \int_\Gamma \lbracket v \rbracket \left \{ k \nabla u \cdot \bm{n} \right \} d\Gamma - \int_\Gamma \left \{ k \nabla v \cdot \bm{n} \right \} \lbracket u \rbracket d\Gamma \nonumber\\
& \quad + \gamma_N \int_\Gamma \lbracket v \rbracket \lbracket u \rbracket d\Gamma =0 \quad \forall v \in V_0 \text{ ,}
\label{eq:nit}
\end{align}
\noindent where $\gamma_N$ is a constraint factor for Nitsche's method.

\subsection{Discretization}

The level set function is discretized by the finite element mesh, such that
\beq
\phi(\bm{x}) = \sum_{i \in I} N_i(\bm{x}) \phi_i \text{ ,}
\eeq
\noindent where $\phi_i$ is the value of the level set function at node $i$. In this work, the interface position is prescribed by determining $\phi_i$ at each node using the signed distance function (\ref{eq:signdist}). Since $\phi(\bm{x})$ is discretized by the finite element mesh, the resolution of the inclusion geometry is dependent on $\mathcal{T}_h$ and improves with mesh refinement. The intersection of $\Gamma$ with an element edge is identified by a sign change in $\phi_i$ for a pair of edge nodes. The intersection of $\Gamma$ directly through a node or an element edge is avoided by enforcing $\phi_i \ne 0$. For any node $i$ where $\left\lvert \phi_i \right\rvert<\phi_{min}$, the nodal level set value is changed to $\phi_i=-\phi_{min}$. For the examples in this work, $\phi_{min}=2 \cdot 10^{-9} \sqrt{\frac{A^e}{\pi}}$ where $A^e$ is the element area. 

\begin{figure}[ht]
\subfloat[]{\includegraphics[trim=1.2in 5in 1.2in 1.8in,clip,width=0.45\textwidth]{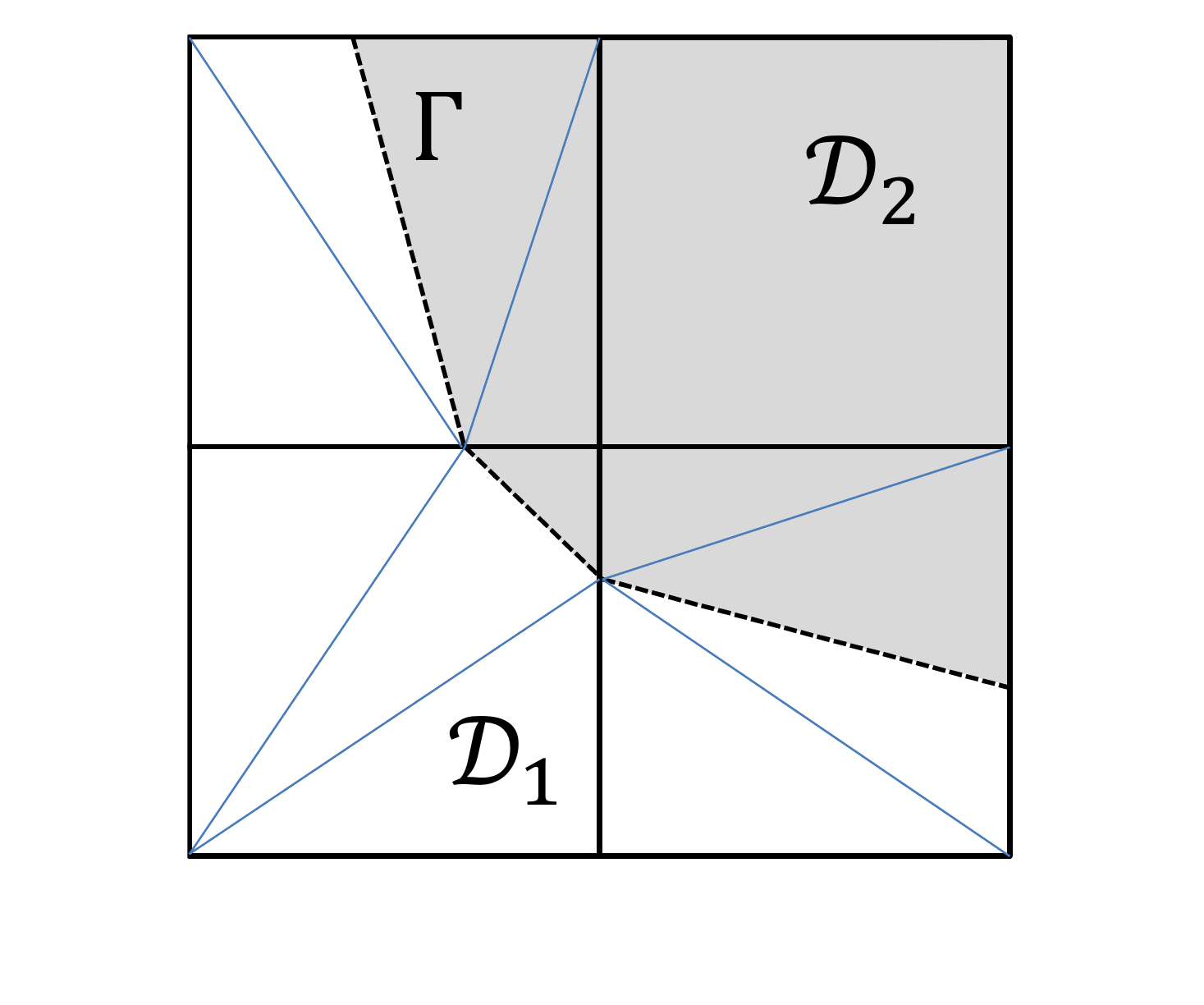}}
\subfloat[]{\includegraphics[trim=1.2in 5in 1.2in 1.8in,clip,width=0.45\textwidth]{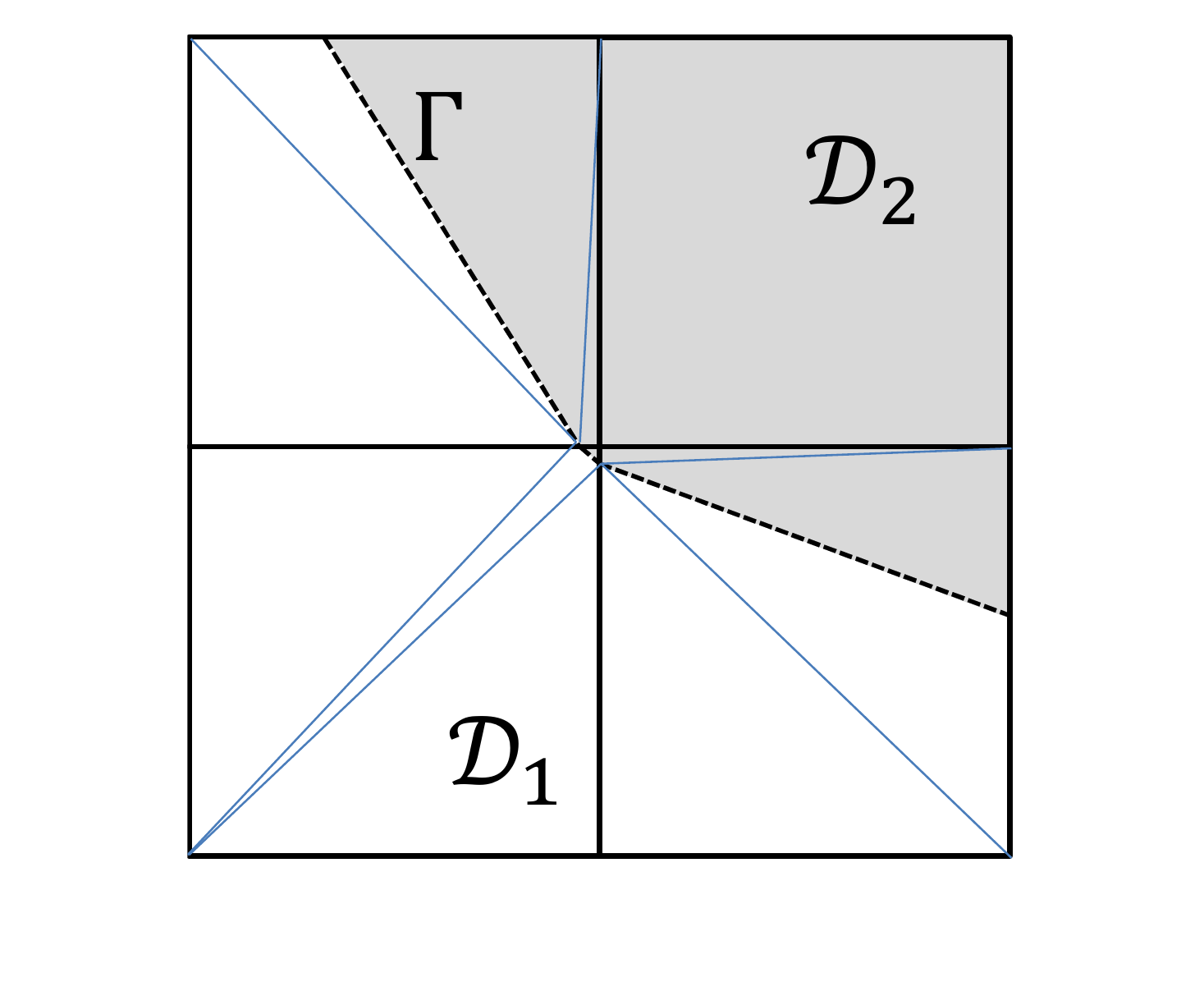}}
\caption{Triangulated partition of a four element configuration leading to a (a) well-conditioned and (b) ill-conditioned system.}
\label{fig:xfemelem}
\end{figure}

Accurate integration over intersected elements is performed by partitioning the element domain, $\mathcal{D}^e$, for piecewise integration. In particular, we partition $\mathcal{D}^e$ using a triangulation aligned with $\Gamma$. An illustration of the triangulation is shown in Fig. \ref{fig:xfemelem} for two configurations of the interface using four elements.

We consider a uniform mesh for $\mathcal{T}_h$ constructed with quadrilateral elements. Bilinear nodal basis functions are used for $N_i(\bm{x})$. For the model problem, elemental Lagrange multipliers are introduced for the stabilized Lagrange multiplier method. We choose a constant approximation of $\lambda$ along the interface $\Gamma$ in an intersected element. This approximation for $\lambda$ is chosen for convenience, as it allows condensing the Lagrange multiplier degree of freedom at an elemental level; other approximations of the Lagrange multiplier can be used in the formulation (\ref{eq:slag}). The third numerical example in Section \ref{sec:examples} approximates the elemental Lagrange multipliers by bilinear shape functions.

The system of equations is constructed by substituting the approximation (\ref{eq:genxfem}) into the weak form (\ref{eq:slag}) or (\ref{eq:nit}). The test functions for the model problem are defined as $v=N_i(\bm{x})$ and $\mu=1$ following the Bubnov-Galerkin method. The integration in (\ref{eq:slag}) or (\ref{eq:nit}) is performed over each element and assembled to construct the system of equations. The discretized system of equations is given by 
\beq
\bm{K}\hat{\bm{u}}=\bm{f} \text{ ,}
\label{eq:kuf}
\eeq
\noindent where $\hat{\bm{u}}$ is the solution vector collecting the degrees of freedom $u_{i,m}^{(1)}$ and $u_{i,m}^{(2)}$, and $\bm{K}$ and $\bm{f}$ are the conduction matrix and load vector, respectively. For the case in which (\ref{eq:kuf}) may be linear or nonlinear, the system residual and Jacobian may be used with the Newton-Raphson method to solve the system. For the remainder of this paper, we refer to the system residual, $\bm{R}$, and Jacobian, $\bm{J}$, defined as
\begin{align}
& \bm{R} = \bm{K}\hat{\bm{u}} - \bm{f} \\
& J_{ij} = \frac{\partial R_i}{\partial u_j} \text{ .}
\end{align}

\noindent Note that for a linear system of equations, $\bm{J}=\bm{K}$ and only one iteration in the Newton-Raphson method is required.

\section{Preconditioning Scheme} \label{sec:precon}

We propose a preconditioning scheme in order to transform the system of equations into a form that is well-conditioned and suitable for solving iteratively. For the configuration shown in Fig. \ref{fig:xfemelem} (b), the lower left element has a small ratio of intersected areas. The region of influence for the degree of freedom at the lower left node interpolating phase 2 approaches zero as the interface approaches the center node. The region of influence for a degree of freedom is the intersection of the nodal basis function support with the physical subdomain $\mathcal{D}_i$. Our aim is to mitigate the sensitivity of the residual to the dissimilar regions of influence for the degrees of freedom. The proposed approach consists of transforming the degrees of freedom by a preconditioning matrix and constraining degrees of freedom associated with small intersections to zero. The constrained degrees of freedom are removed from the equations when solving the system. When constraining degrees of freedom only without using the preconditioning matrix \citep{ref:reusken2008}, the solution accuracy decreases as the condition number is reduced. The proposed approach solves the problem in a transformed space and does not change the solution to the discrete problem. We will show in Section \ref{sec:examples} that the preconditioning scheme maintains an approximately constant condition number without loss of solution accuracy.

A geometric preconditioner $\bm{T}$ is introduced, such that the solution in the physical space, $\hat{\bm{u}}$, is obtained by

\beq
\hat{\bm{u}}=\bm{T}\tilde{\bm{u}} \text{ ,}
\label{eq:solproj}
\eeq

\noindent where $\tilde{\bm{u}}$ is the solution in the transformed space. The residual and Jacobian of the system in the transformed space are defined as
\begin{align}
\tilde{\bm{R}} & = \bm{T}^T \bm{R} \nonumber \\
\tilde{\bm{J}} & = \bm{T}^T \bm{J} \bm{T} \text{ .}
\label{eq:proj}
\end{align}

Note, the residual $\bm{R}$ and the Jacobian $\bm{J}$ are constructed in a standard fashion using the XFEM. For problems with dynamically evolving interfaces, such as phase change and multi-phase flow problems \citep{ref:chessa2002,ref:zabaras2006,ref:chessa2003}, the discretized level set field contributes degrees of freedom to the solution vectors $\hat{\bm{u}}$ and $\tilde{\bm{u}}$. In this case, the Jacobian $\tilde{\bm{J}}$ contains additional terms. We omit a detailed discussion of this class of problems and focus on problems with static or prescribed interface geometries.

The purpose of the geometric preconditioner is to balance the influence for degrees of freedom as the intersected areas approach zero. There are two issues associated with the intersected areas approaching zero. First, the partitioned element integration, and therefore the diagonal entry of the element matrix, approaches zero because the area of integration is small. Second, the influence of a degree of freedom on the residual will vanish as the region of influence approaches zero.

Here, we construct a diagonal preconditioning matrix for $\bm{T}$ from the nodal basis functions and their support in order to transform the degrees of freedom. The proposed approach accommodates other choices for $\bm{T}$, both diagonal and non-diagonal. However, diagonal scaling is more computationally efficient in terms of memory and matrix operations. The preconditioning matrix $\bm{T}$ is constructed by integrating the nodal basis functions ($\bm{T}_N$) or derivatives ($\bm{T}_B$) over the nodal support. The diagonal components of the $\bm{T}_N$ preconditioning matrix are defined as
\beq
T_{i,m}^{(p)} = \left( \max_{e \in E_i} \frac{\int_{\mathcal{D}_p^e} N_i(\bm{x}) dx}{\int_{\mathcal{D}^e} N_i(\bm{x}) dx} \right)^{-\frac{1}{2}} \text{ ,}
\label{eq:ntypemax}
\eeq
\noindent where $T_{i,m}^{(p)}$ corresponds to the degree of freedom $u_{i,m}^{(p)}$ at node $i$, and $E_i$ is the set of elements connected to node $i$. Here, $\mathcal{D}_p^e$ denotes the element domain which belongs to phase $p$. The diagonal components of the $\bm{T}_B$ preconditioning matrix are defined as
\beq
T_{i,m}^{(p)} = \left( \max_{e \in E_i} \frac{\int_{\mathcal{D}_p^e} \nabla N_i(\bm{x}) \cdot \nabla N_i(\bm{x}) dx}{\int_{\mathcal{D}^e} \nabla N_i(\bm{x}) \cdot \nabla N_i(\bm{x}) dx} \right)^{-\frac{1}{2}} \text{ .}
\label{eq:btypemax}
\eeq

\noindent In practice, the components $T_{i,m}^{(p)}$ are only computed at nodes connected to an intersected element. If all elements in $E_i$ are non-intersected, then the degrees of freedom at node $i$ are not transformed and $T_{i,m}^{(p)}=1$.

The $\bm{T}_N$ and $\bm{T}_B$ geometric preconditioners both lead to scaling terms that increase as the region of influence for degrees of freedom approaches zero. The region of influence is measured by $\max_{e \in E_i}{\int_{\mathcal{D}_p^e} N_i(\bm{x}) dx}$ and $\max_{e \in E_i}{\int_{\mathcal{D}_p^e} \nabla N_i(\bm{x}) \cdot \nabla N_i(\bm{x}) dx}$ in (\ref{eq:ntypemax}) and (\ref{eq:btypemax}), respectively. For a given problem, the choice of the preconditioner type can be determined by the dominating operator in the partial differential equation. Based on the construction of the system of equations, the $\bm{T}_B$ preconditioning matrix is more appropriate for diffusion dominated problems, while $\bm{T}_N$ is appropriate for convection or reaction dominated problems.

The preconditioner $\bm{T}$ improves the condition number by balancing the influence of the degrees of freedom. However, as the preconditioner is constructed using the nodal basis functions, the scaling terms in $\bm{T}$ do not approach $\infty$ at the same rate as the region of influence approaches zero. Therefore, an ill-conditioned system of equations may still result when the ratio of intersected areas approaches zero. In addition to the preconditioner, we propose to constrain degrees of freedom to zero with small regions of influence. The criteria for selecting the degrees of freedom to be constrained to zero is defined as
\beq
T_{i,m}^{(p)} > T_{tol} \text{ ,}
\label{eq:Ttol}
\eeq

\noindent where $T_{tol}$ is a specified tolerance. It is shown in Section \ref{sec:examples} that there is a wide range for the choice of $T_{tol}$ which does not impact the numerical error and condition number. Constraining degrees of freedom to zero is needed when $T_{i,m}^{(p)}>>1$. The numerical studies in Section \ref{sec:examples} suggest values for $T_{tol}$ between $10^{4}$ and $10^{8}$.

A summary of applying the proposed preconditioning scheme to a nonlinear problem solved by the Newton-Raphson method is outlined below:

\begin{enumerate}
  \item Construct $\mathcal{T}_h$ and $\phi$.
  \item Construct $\bm{T}$ using (\ref{eq:ntypemax}) or (\ref{eq:btypemax}) and mark degrees of freedom to be constrained by (\ref{eq:Ttol}).
  \item Obtain transformed initial guess by the inverse operation of (\ref{eq:solproj}).
  \item Solve iteratively the problem $\tilde{\bm{R}}=0$ for $\tilde{\bm{u}}$ as follows:
  \begin{enumerate}
  \item Reconstruct $\bm{T}$ and update degrees of freedom to be constrained.
    \item Obtain current solution by (\ref{eq:solproj}).
    \item Construct $\bm{R}$ and $\bm{J}$.
    \item Obtain $\tilde{\bm{R}}$ and $\tilde{\bm{J}}$ by (\ref{eq:proj}).
    \item Solve transformed system for $\Delta\tilde{\bm{u}}$.
		\item Update solution and check for convergence.
  \end{enumerate}
  \item Obtain final solution, $\hat{\bm{u}}$, by (\ref{eq:solproj})
\end{enumerate}

As shown in the implementation outline, $\bm{T}$ is constructed prior to computing the residual and Jacobian. If the interface geometry is prescribed and independent of the solution, then the level set field and hence $\bm{T}$ do not change in the Netwon iterations. In this case, step 4(a) is not necessary.

\section{Numerical Examples}
\label{sec:examples}

In this section, the performance of the preconditioning scheme is studied for three problems. The first example illustrates the basic concept of the preconditioning scheme when solving a diffusion problem for a two-material bar. The second example is a diffusion problem with a circular material inclusion. For these examples, the accuracy of the solution as well as the condition number of the systems are examined with and without the proposed preconditioning scheme. The third example is a transient flow problem with a moving rigid obstacle, modeled by the incompressible Navier-Stokes equations. This example demonstrates the applicability of the proposed scheme to nonlinear transient problems with moving interfaces. While the examples in this paper consider 2D problems, the extension of the proposed preconditioning scheme to 3D problems is straight-\linebreak forward.

\subsection{Example 1: Two-Material Bar Diffusion}

\begin{figure}[ht]
\centering\includegraphics[trim=1.5in 6.4in 0.4in 2.1in,clip,width=0.6\textwidth]{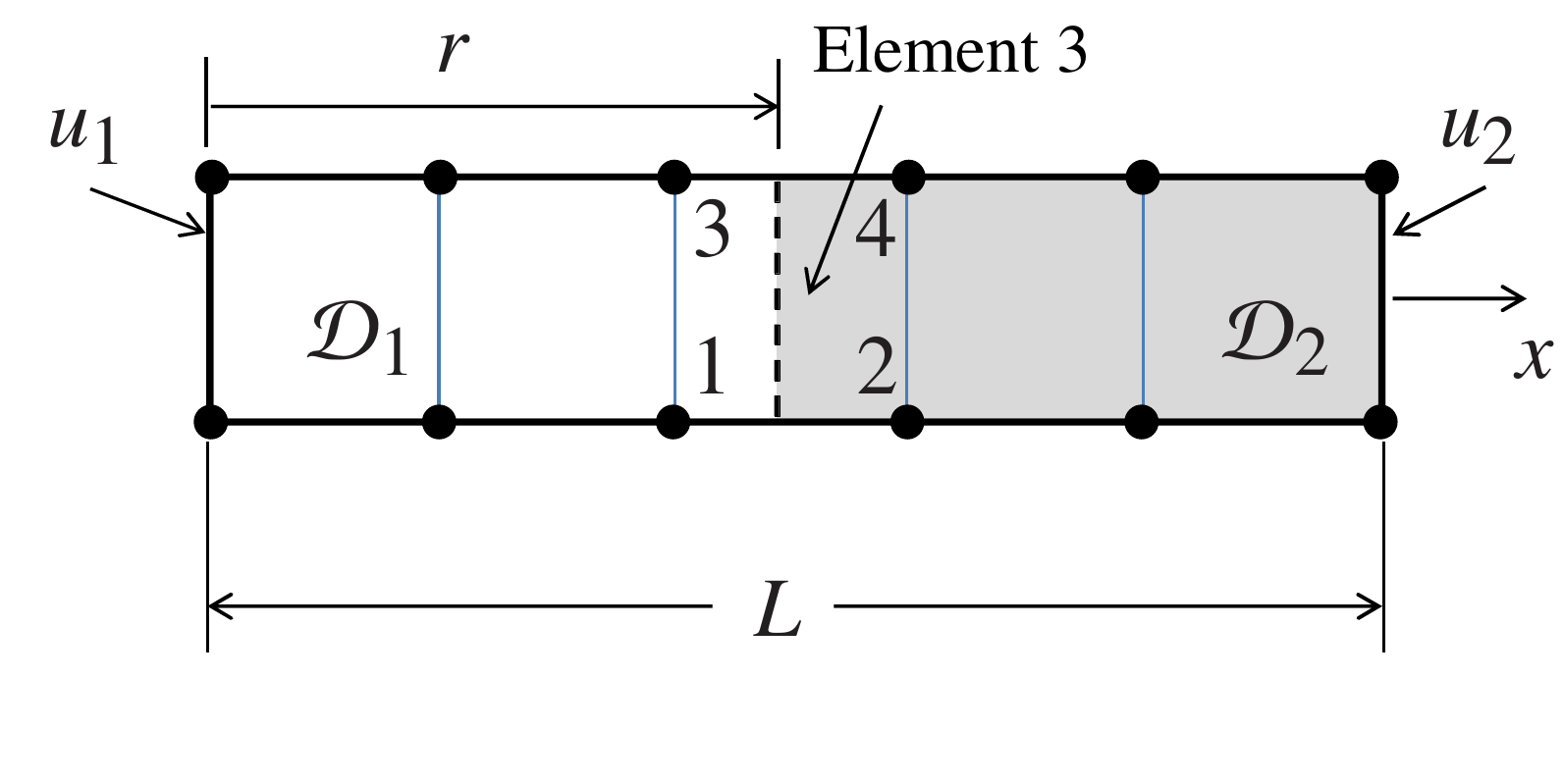}
\caption{Problem description for Example 1.}
\label{fig:bar}
\end{figure}

We illustrate the basic concept of the preconditioning scheme for a simple example with an analytical solution. We consider solving the heat conduction model for the two-material bar shown in Fig. \ref{fig:bar}. The length of the bar is $L$, and temperatures $u_1$ and $u_2$ are specified at $x=0$ and $x=L$, respectively. The material conductivity is $k_1=1$ in $\mathcal{D}_1$ and $k_2=2$ in $\mathcal{D}_2$. The position of the vertical interface is measured from the left end and specified by $r$. The problem is solved using quadrilateral elements. While the exact solution can be captured using one element, we discretize the bar with five elements in order to vary the position of the interface across one element. Note that while this example is useful for explaining the concept and demonstrating the reduced condition number, it is not well suited to illustrate a change in the accuracy of the solution due to an ill-conditioned system. Without preconditioning, an ill-conditioned system will occur when the interface is nearly aligned with an element edge. In the intersected element, a ratio of the area of the phase 1 and phase 2 regions with a value less than $10^{-13}$ results in a condition number greater than $10^{14}$.

The interface position is varied from $r/L=0.3$ to $r/L=0.7$ in steps of $\Delta r/L=0.002$. Element 3 is intersected for $0.4<r/L<0.6$. As $r/L$ approaches $0.4$ and $0.6$, the ratio of intersected areas in element 3 becomes small. The preconditioning scheme using the $\bm{T}_B$ matrix and $T_{tol}=10^4$ is selected for the example bar problem using the stabilized Lagrange multiplier method with $\gamma_S=k_1+k_2$ for enforcing continuity at the interface.

There are four degrees of freedom for element 3 at nodes 1 to 4 which have small regions of influence as the interface position is varied. Since the problem is one-dimensional, we only consider nodes 1 and 2, and focus on the degrees of freedom $u_{1,1}^{(2)}$ and $u_{2,1}^{(1)}$. The degree of freedom $u_{1,1}^{(2)}$ is used for interpolating the phase 2 solution in element 3, and it has a small region of influence when $r/L \approx 0.6$. The degree of freedom $u_{2,1}^{(1)}$ is used for interpolating the phase 1 solution in element 3, and it has a small region of influence when $r/L \approx 0.4$.

\begin{figure}[ht]
\centering
\subfloat[]{\includegraphics[trim=1.6in 5in 1in 2.05in,clip,width=0.6\textwidth]{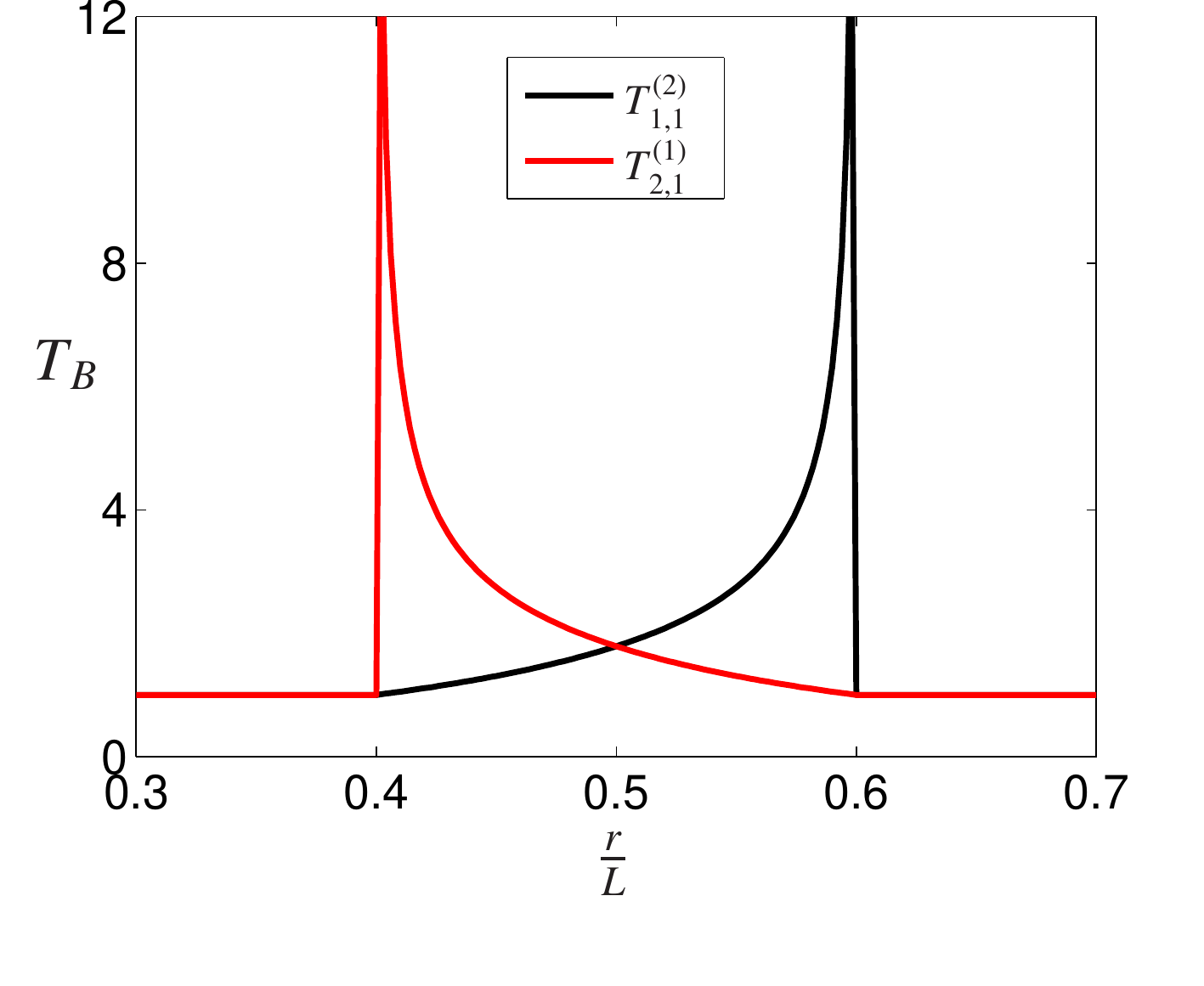}} \\
\subfloat[]{\includegraphics[trim=1.6in 5in 1in 2.05in,clip,width=0.6\textwidth]{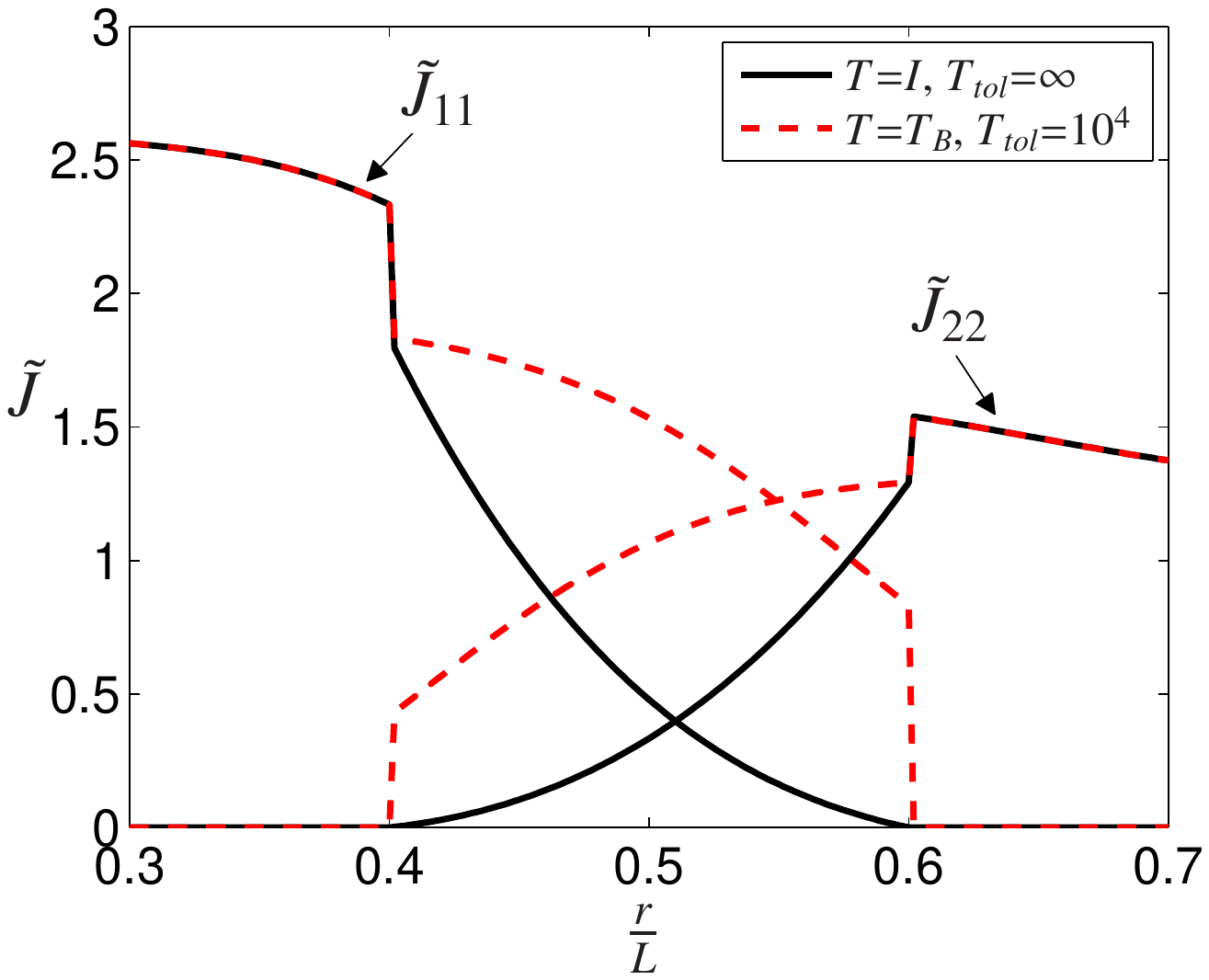}}
\caption{The diagonal components of (a) $\bm{T}_B$ and (b) $\tilde{\bm{J}}$ corresponding to the degrees of freedom $u_{1,1}^{(2)}$ and $u_{2,1}^{(1)}$. Here $\tilde{J}_{11}$ and $\tilde{J}_{22}$ correspond to $u_{1,1}^{(2)}$ and $u_{2,1}^{(1)}$, respectively.}
\label{fig:Tmatbar}
\end{figure}

The $\bm{T}_B$ values corresponding to these degrees of freedom are shown in Fig. \ref{fig:Tmatbar}(a) as the interface location varies. The $T_{1,1}^{(2)}$ and $T_{2,1}^{(1)}$ values increase as the ratio of intersected areas in element 3 decrease. The diagonal components of $\tilde{\bm{J}}$ corresponding to $u_{1,1}^{(2)}$ and $u_{2,1}^{(1)}$ without preconditioning ($\bm{T}=\bm{I}$) and with the preconditioner $\bm{T}_B$ are shown in Fig. \ref{fig:Tmatbar}(b). The diagonal components of $\tilde{\bm{J}}$ with $\bm{T}=\bm{T}_B$ do not reduce to zero as the ratio of intersected areas approach zero. The jumps in $\tilde{J}_{ii}$ at $r/L=0.4$ and $r/L=0.6$ result from the stabilized Lagrange method for enforcing continuity at the interface.

\begin{figure}[ht]
\centering\includegraphics[trim=1.5in 5.0in 0.3in 2.05in,clip,width=0.65\textwidth]{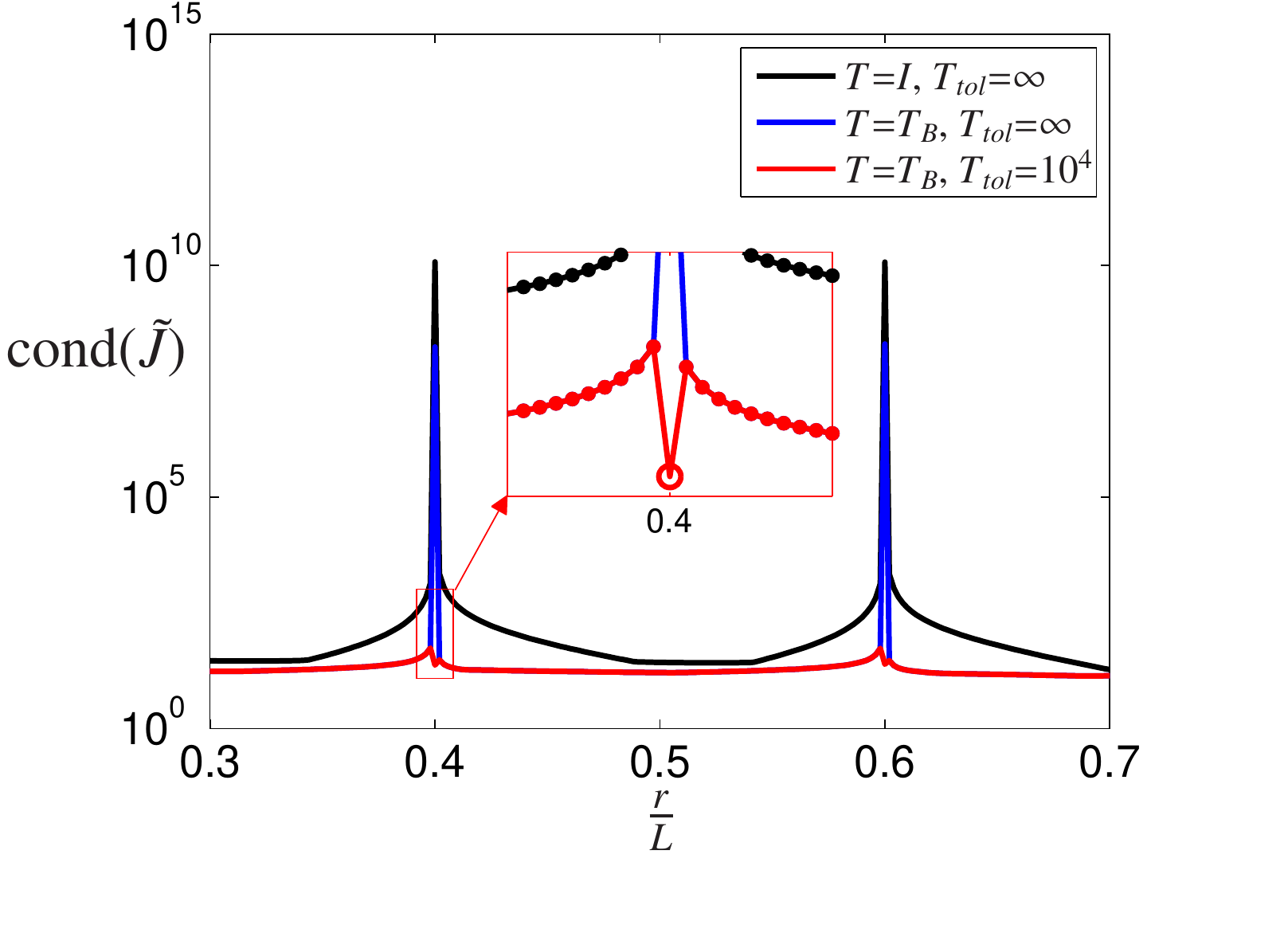}
\caption{Condition number as a function of the interface position for the two-material bar. In the inset figure, $\Delta r/L=2\cdot10^{-5}$ and the open circle marks the $r/L$ value for which degrees of freedom were constrained.}
\label{fig:condbar}
\end{figure}

\begin{figure}[ht]
\centering
\subfloat[]{\includegraphics[trim=1.6in 5in 0.7in 2.05in,clip,width=0.6\textwidth]{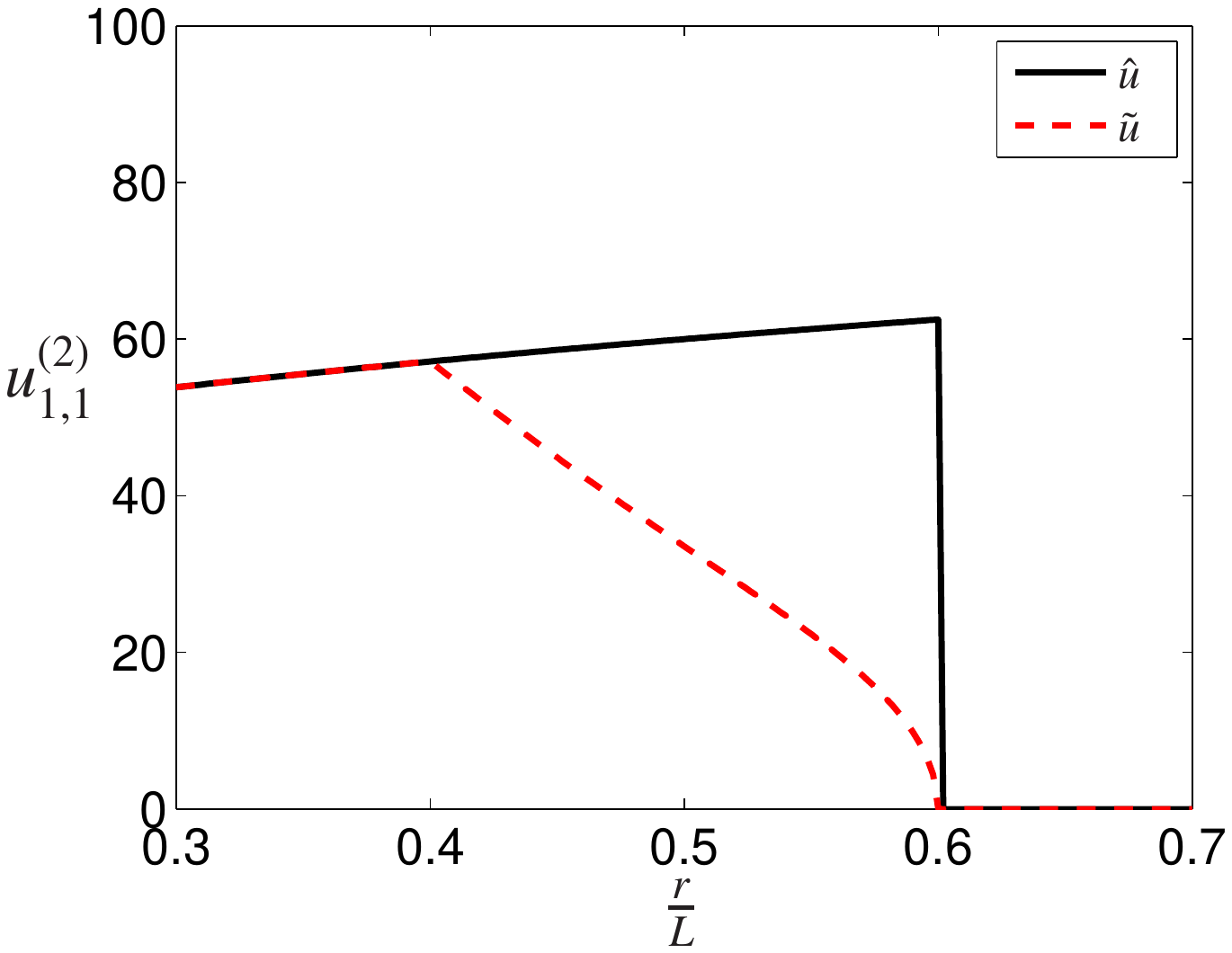}} \\
\subfloat[]{\includegraphics[trim=1.6in 5in 0.7in 2.05in,clip,width=0.6\textwidth]{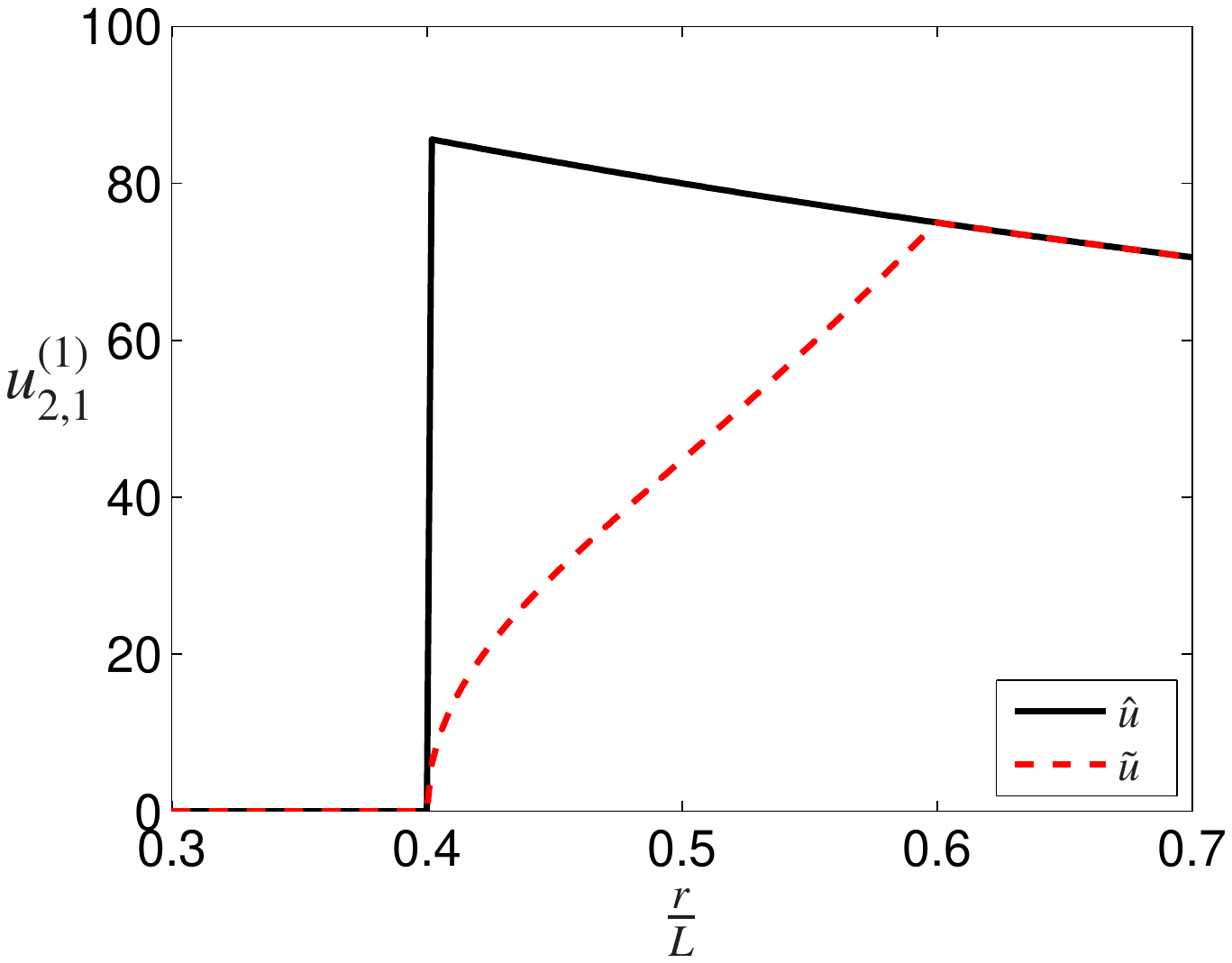}}
\caption{The physical ($\hat{u}$) and projected ($\tilde{u}$) solutions for the degrees of freedom (a) $u_{1,1}^{(2)}$ and (b) $u_{2,1}^{(1)}$.}
\label{fig:bardofs}
\end{figure}

The condition number of $\tilde{\bm{J}}$ is shown as a function of the interface position in Fig. \ref{fig:condbar}. The condition number was determined without and with the preconditioning scheme, denoted by $\bm{T}=\bm{I}$ and $\bm{T}=\bm{T}_B$, respectively. No degrees of freedom were constrained for $T_{tol}=\infty$. The condition number is improved for $\bm{T}=\bm{T}_B$ and $T_{tol}=\infty$, but is still large near $r/L=0.4$ and $r/L=0.6$. By imposing the criteria for constraining degrees of freedom, the condition number at $r/L=0.4$ and $r/L=0.6$ is significantly reduced. The physical and transformed solutions for the degrees of freedom $u_{1,1}^{(2)}$ and $u_{2,1}^{(1)}$ are shown in Fig. \ref{fig:bardofs}. The physical degrees of freedom jump to zero when element 3 is not intersected. The influence of the preconditioning for $u_{1,1}^{(2)}$ and $u_{2,1}^{(1)}$ occurs when element 3 is intersected. The projected degrees of freedom vary to zero as the ratio of intersected areas approach zero.

\subsection{Example 2: Circular Inclusion Diffusion}

\begin{figure}[ht]
\centering\includegraphics[trim=3in 1.5in 3in 1.5in,clip,width=0.55\textwidth]{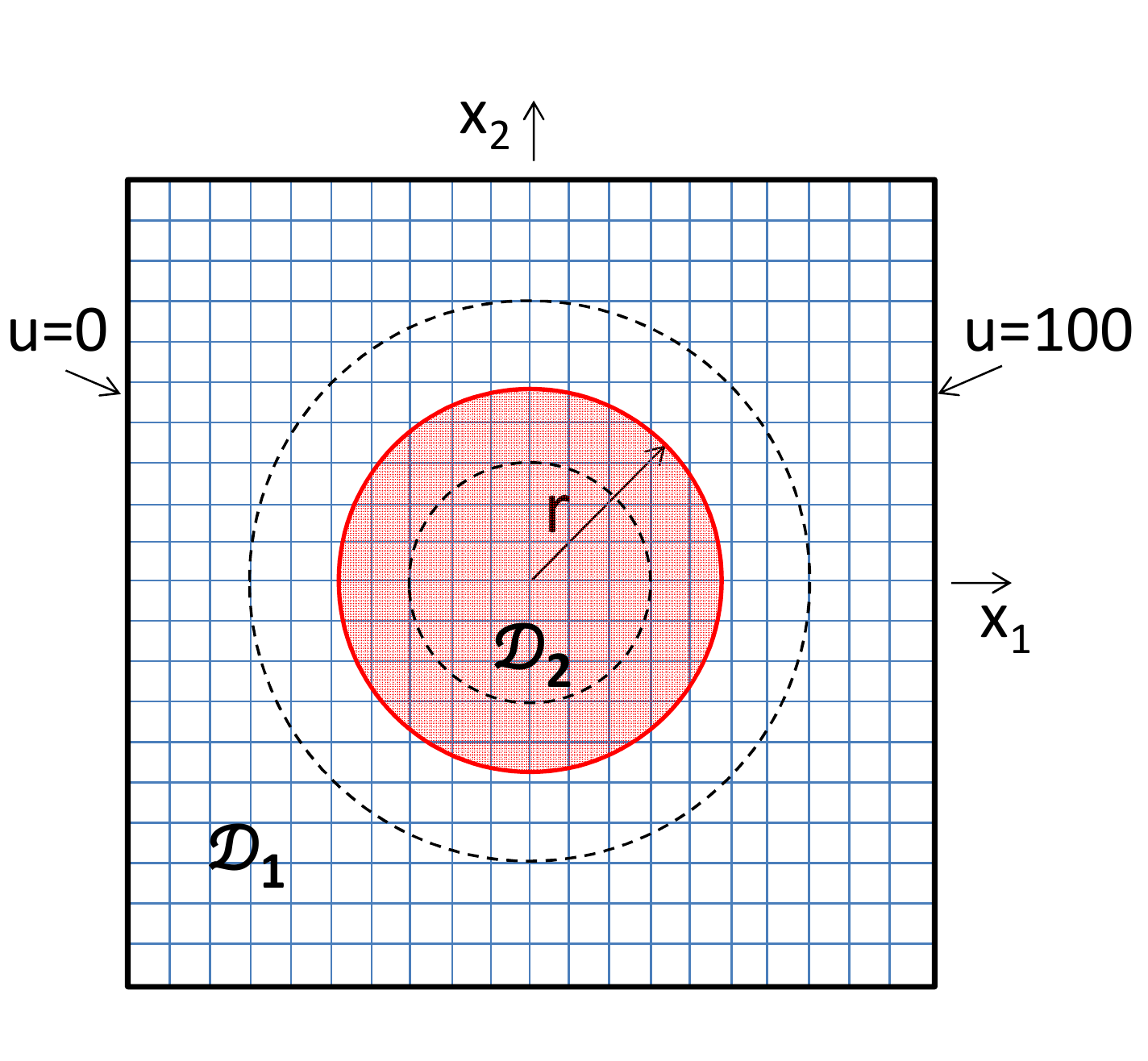}
\caption{Problem description for Example 2.}
\label{fig:ex1}
\end{figure}

The second numerical example is the heat transfer problem shown in Fig. \ref{fig:ex1}. The model problem (\ref{eq:model}) is solved for a square domain $\mathcal{D}=(-10,10)\times(-10,10)$ with a centered circular inclusion of radius $r$. The radius is varied from $r=3$ to $r=7$ in steps of $\Delta r=0.02$. Material 1 has a conductivity $k_1=2$ in $\mathcal{D}_1$, and material 2 has a conductivity $k_2=2\cdot10^3$ in $\mathcal{D}_2$. The temperature is specified as $u=0$ on the left boundary and $u=100$ on the right boundary. The top and bottom edges are adiabatic. The two methods of enforcing the solution continuity at the interface (\ref{eq:slag}) and (\ref{eq:nit}) are considered with $\gamma_S=k_1+k_2$ and $\gamma_N=10^{-3}(k_1+k_2)$.

The condition number of the system of equations depends on the configuration of the intersections and the ratio of conductivities. A high ratio of conductivities, also considered in \citep{ref:mandel1996,ref:ewing2001}, is used here to highlight the ill-conditioning issue for this simple example problem. The ratio of intersected areas is examined for the variation of the radius by determining the minimum element area ratio, defined as
\beq
A_{min}=\min\limits_{e \in \mathcal{T}_h} \frac{\mathcal{D}_1^e}{\mathcal{D}_2^e} \text{ .}
\eeq

\begin{figure}[ht]
\centering\includegraphics[trim=1.5in 5in 0.8in 2in,clip,width=0.55\textwidth]{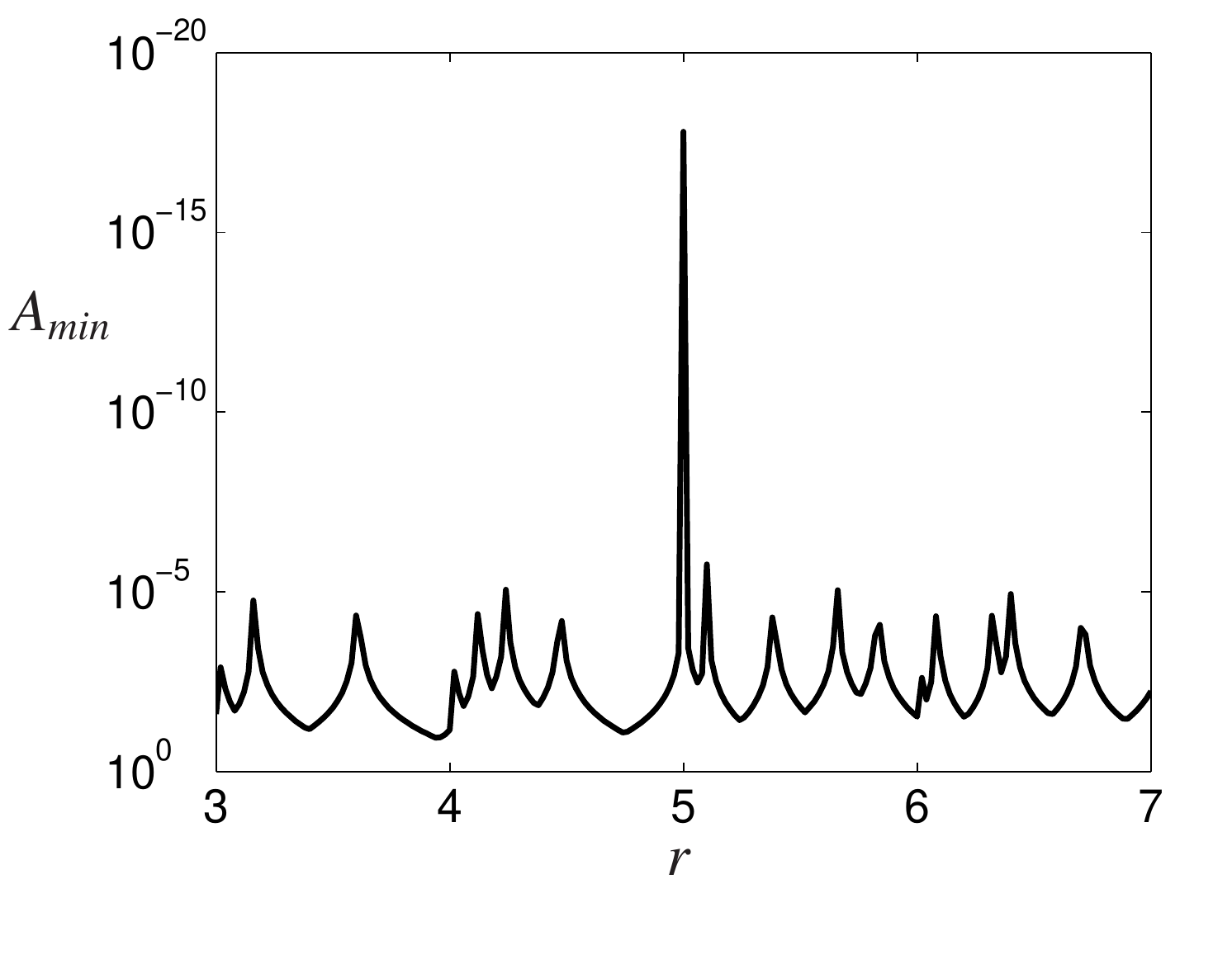}
\caption{Minimum element area ratio, $A_{min}$, for Example 2.}
\label{fig:minarea}
\end{figure}

\noindent The variation of $A_{min}$ with $r$ is shown in Fig. \ref{fig:minarea}. Note the vertical axis is reversed, such that small intersections are indicated by the peaks. The minimum area ratios of order $10^{-5}$ and $10^{-18}$ which occur for the variation of $r$ lead to a high condition number of the system.

Three studies were performed for this example. The first study shows the influence of $T_{tol}$ in (\ref{eq:Ttol}) on the condition number and solution accuracy. The second study is a comparison of the condition number using a body-fitted mesh, XFEM with a Jacobi preconditioner, and XFEM with the proposed preconditioning scheme. Finally, we study the influence of the preconditioning scheme on the performance of an iterative solver.

To study the influence of $T_{tol}$ on the condition number and solution accuracy, the $\bm{T}_B$ preconditioning matrix and the stabilized Lagrange multiplier method are used. The value of $T_{tol}$ is varied from $T_{tol}=10$ to $T_{tol}=10^8$. The maximum condition number of $\tilde{\bm{J}}$ and solution error is computed for each value of $T_{tol}$ by considering all values of $r$. The maximum condition number, $c_{max}$, is defined by 
\beq
c_{max}=\max\limits_r \; \text{cond}(\tilde{\bm{J}}) \text{ .}
\label{eq:maxcond}
\eeq

\begin{figure}[ht]
\centering
\subfloat[]{\includegraphics[trim=1.5in 5in 0.5in 2.05in,clip,width=0.65\textwidth]{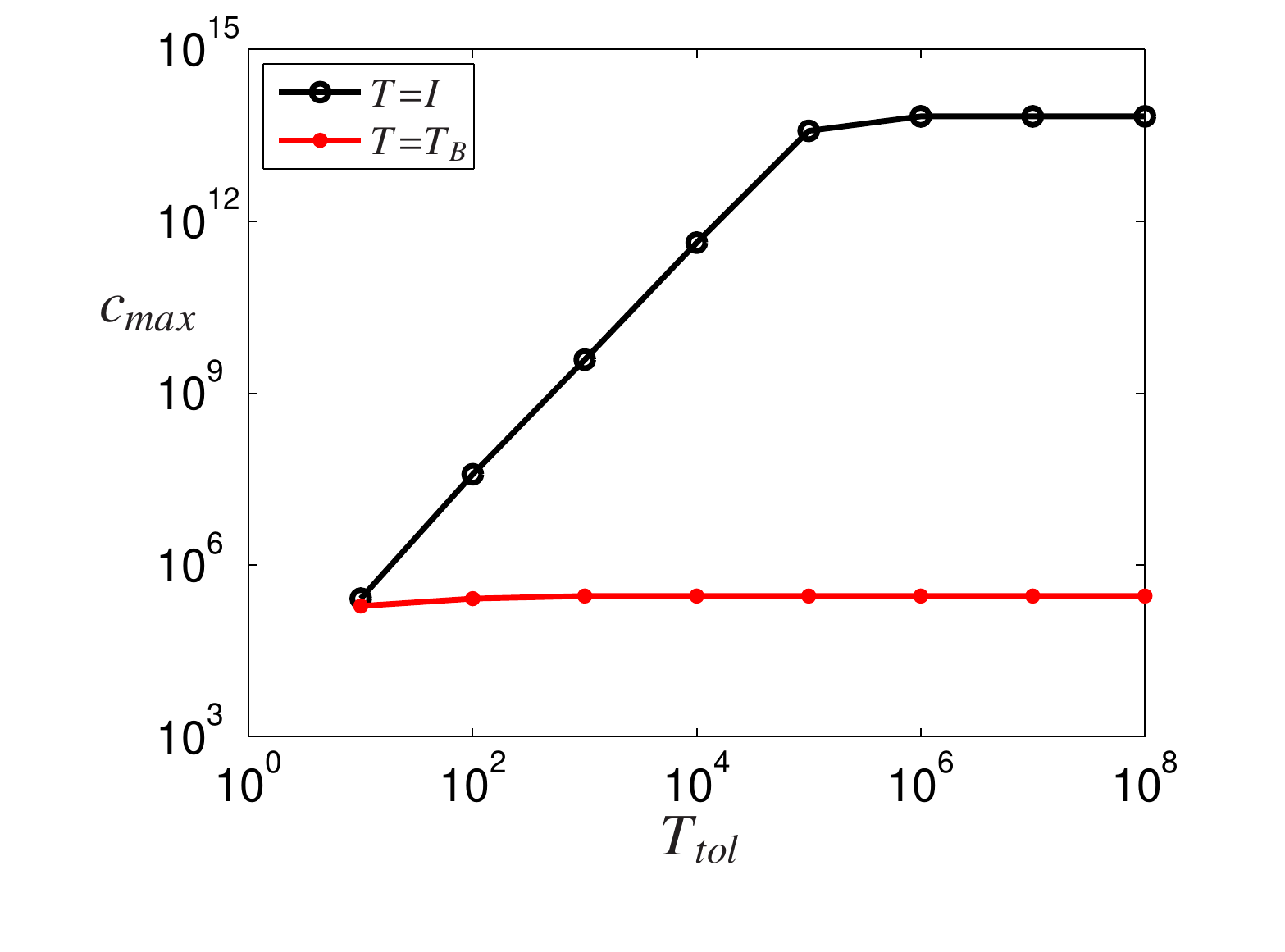}} \\
\subfloat[]{\includegraphics[trim=1.5in 5in 0.5in 2.05in,clip,width=0.65\textwidth]{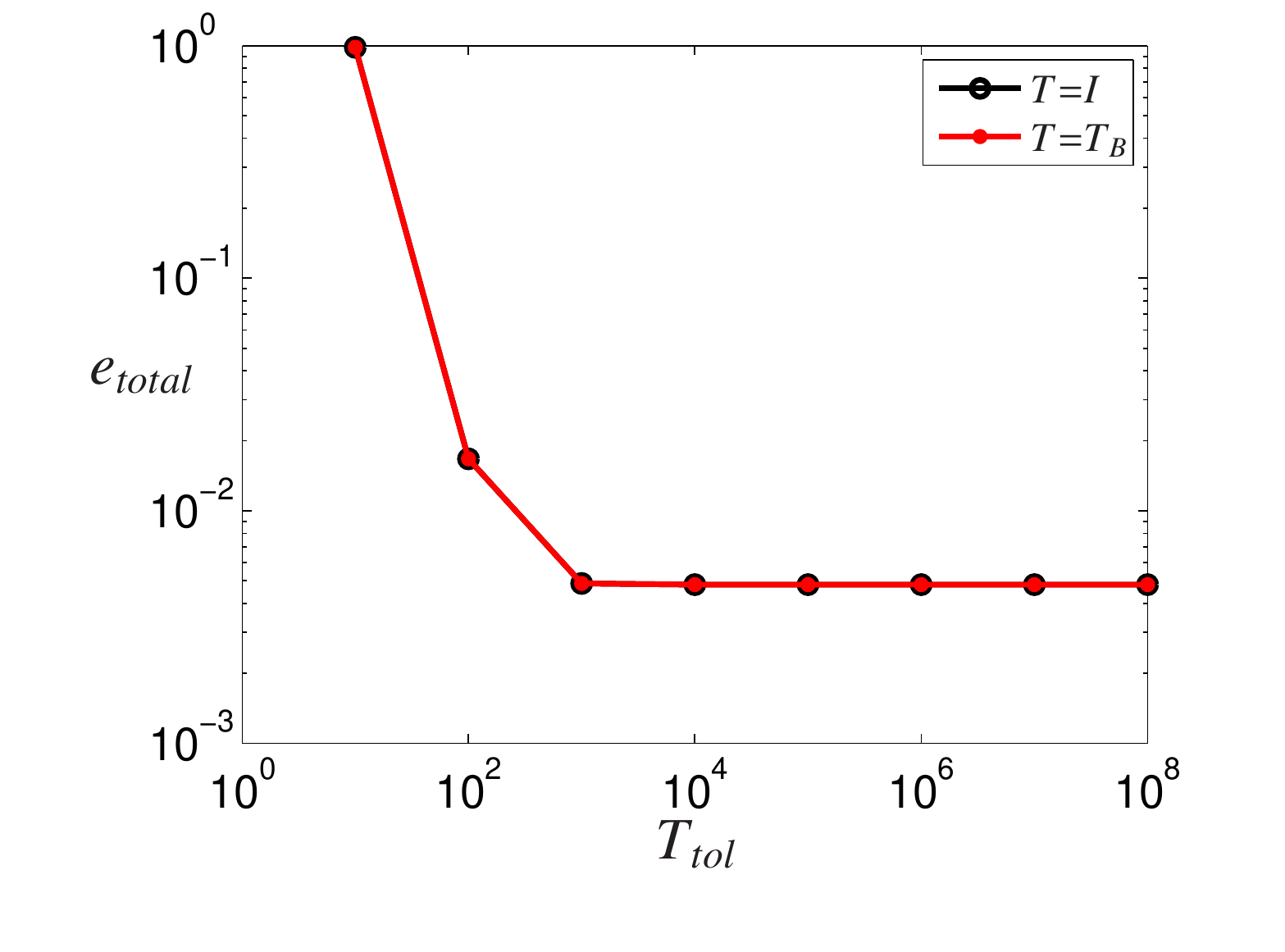}}
\caption{Influence of $T_{tol}$ on the (a) maximum condition number and (b) approximation error computed from (\ref{eq:maxcond})-(\ref{eq:L2err}) for Example 2.}
\label{fig:conderrcir}
\end{figure}

\noindent The accuracy of the XFEM solution is measured by integrating the $L_2$ relative error, such that the total error for each value of $T_{tol}$ is defined by
\beq
e_{total}=\int_3^7{\frac{\| \hat{u}(r)-u_{ref1}(r) \|_2}{\| u_{ref1}(r) \|_2}}dr \text{ ,}
\label{eq:L2err}
\eeq
\noindent where $u_{ref1}(r)$ is a reference solution for radius $r$ obtained using a body-fitted finite element mesh with an element size of $h \approx 0.05$. The influence of $T_{tol}$ on the condition number and solution error is shown in Fig. \ref{fig:conderrcir} with and without the preconditioning matrix. For $\bm{T}=\bm{I}$, the preconditioning matrix is only used for the criteria on constraining degrees of freedom in (\ref{eq:Ttol}) and not applied when solving the system of equations. More degrees of freedom are constrained to zero by decreasing $T_{tol}$, and the maximum condition number is reduced for $\bm{T}=\bm{I}$. However, the solution error increases as more degrees of freedom are constrained. For $\bm{T}=\bm{T}_B$, the condition number is reduced for each value of $T_{tol}$. The solution error is the same for $\bm{T}=\bm{I}$ and $\bm{T}=\bm{T}_B$. Note that the same number of degrees of freedom were constrained to zero for $T_{tol}$ values of $10^6$, $10^7$, and $10^8$.

The second study compares the condition number for various choices of $\bm{T}$. The condition number of $\tilde{\bm{J}}$ is computed for the variation of $r$ using XFEM with the stabilized Lagrange and Nitsche methods. The condition number of $\tilde{\bm{J}}$ using a body-fitted mesh with an element size of $h\approx0.5$ and $\bm{T}=\bm{I}$ was also computed. A Jacobi preconditioner is implemented by defining
\beq
\bm{T}_{jac}=\text{diag}(\bm{J})^{-\frac{1}{2}} \text{ .}
\label{eq:Tjac}
\eeq

\begin{figure}[ht]
\centering\includegraphics[trim=1.4in 5.5in 0.6in 2.0in,clip,width=0.65\textwidth]{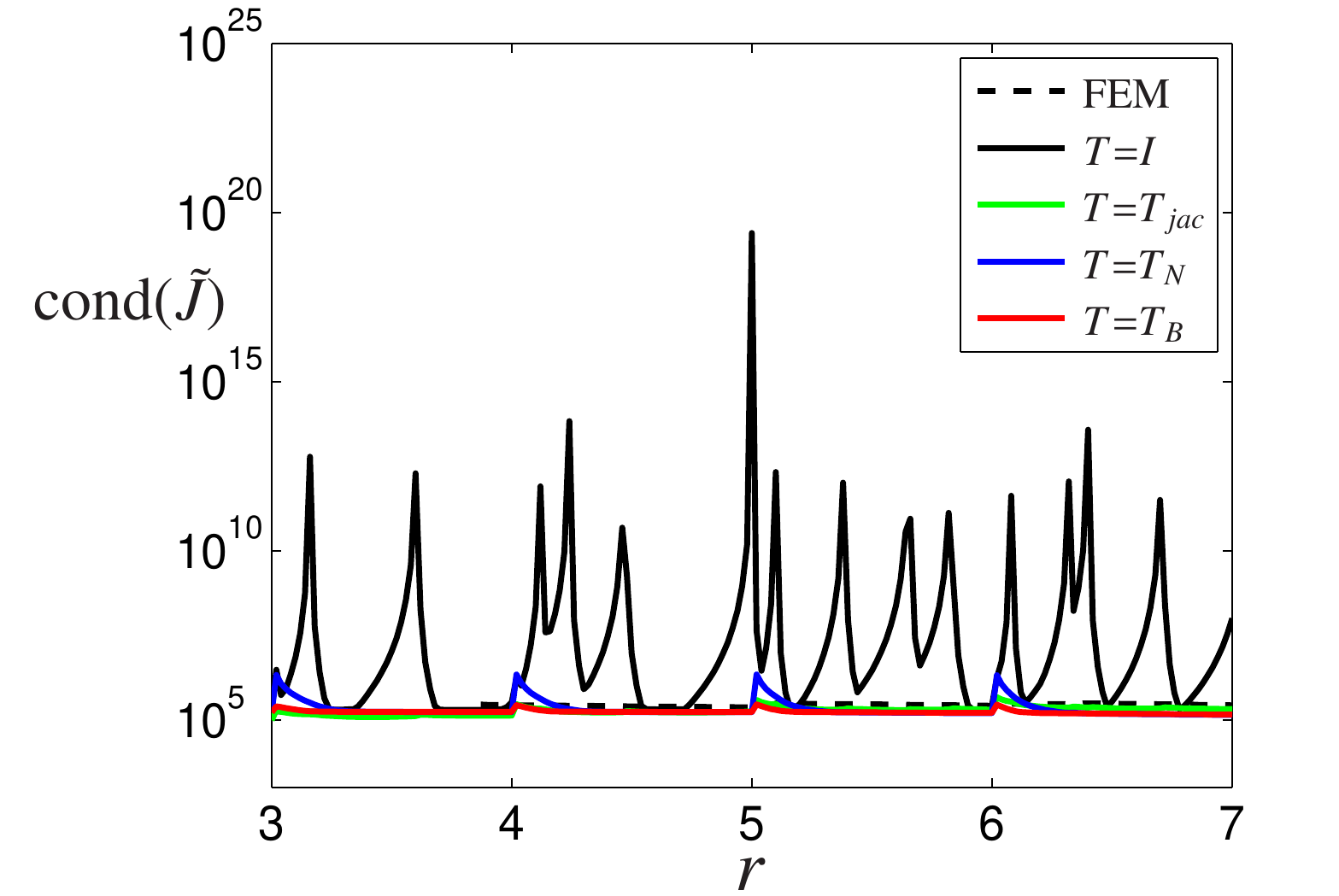}
\caption{Comparison of the condition number for a varying inclusion radius using the stabilized Lagrange method for Example 2.}
\label{fig:condslag}
\end{figure}

\begin{figure}[ht]
\centering\includegraphics[trim=1.4in 5.5in 0.6in 2.0in,clip,width=0.65\textwidth]{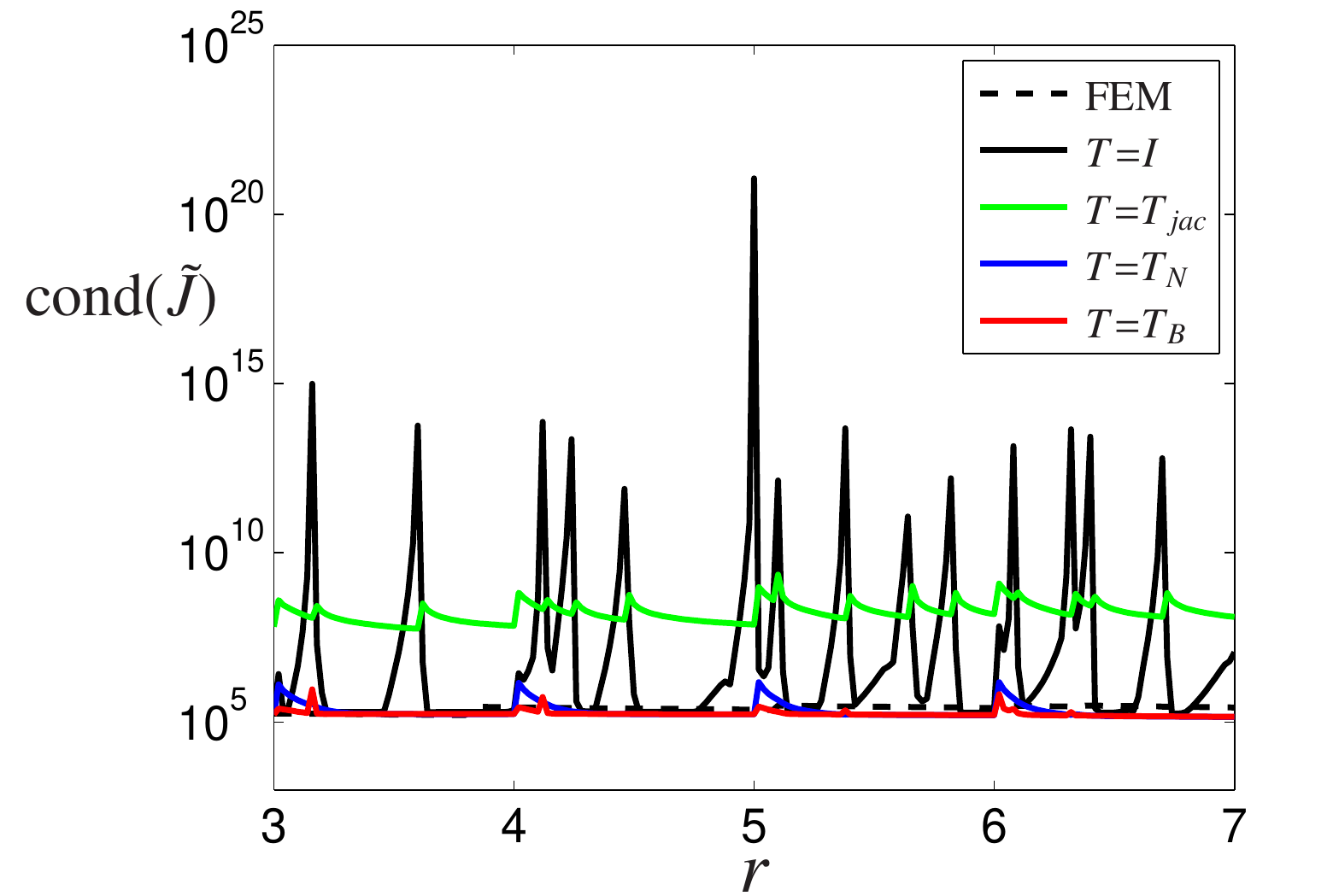}
\caption{Comparison of the condition number for a varying inclusion radius using Nitsche's method for Example 2.}
\label{fig:condnits}
\end{figure}

\noindent Note that $\bm{T}_{jac}$ is a solver preconditioner applicable to solving the linear system and not a geometric preconditioner as used in the proposed scheme. Finally, the condition number of $\tilde{\bm{J}}$ is computed using the $\bm{T}_N$ and $\bm{T}_B$ preconditioning matrices with $T_{tol}=10^8$. A comparison of the condition numbers for the variation of $r$ is shown in Figs. \ref{fig:condslag} and \ref{fig:condnits}. No degrees of freedom were constrained for $\bm{T}=\bm{I}$ and $\bm{T}=\bm{T}_{jac}$, which corresponds to $T_{tol}=\infty$. For $\bm{T}=\bm{I}$, the condition number using XFEM varies with the size of the inclusion up to an order of $10^{20}$. The $r$ values of the high condition numbers correspond to the small intersections seen in Fig. \ref{fig:minarea}. For $\bm{T}=\bm{T}_{jac}$, the condition number is comparable to that of the body-fitted FEM system for the stabilized Lagrange method (Fig. \ref{fig:condslag}) but not Nitsche's method (Fig. \ref{fig:condnits}). This suggests that the condition number is influenced by the off-diagonal terms in $\bm{J}$ for Nitsche's method. However, for $\bm{T}=\bm{T}_N$ and $\bm{T}=\bm{T}_B$, the XFEM condition number is comparable to the body-fitted FEM system for both stabilized Lagrange and Nitsche methods for all interface positions.

\begin{figure*}[ht]
\centering
\subfloat[]{\includegraphics[trim=1.4in 5.2in 1in 2.1in,clip,width=0.33\textwidth]{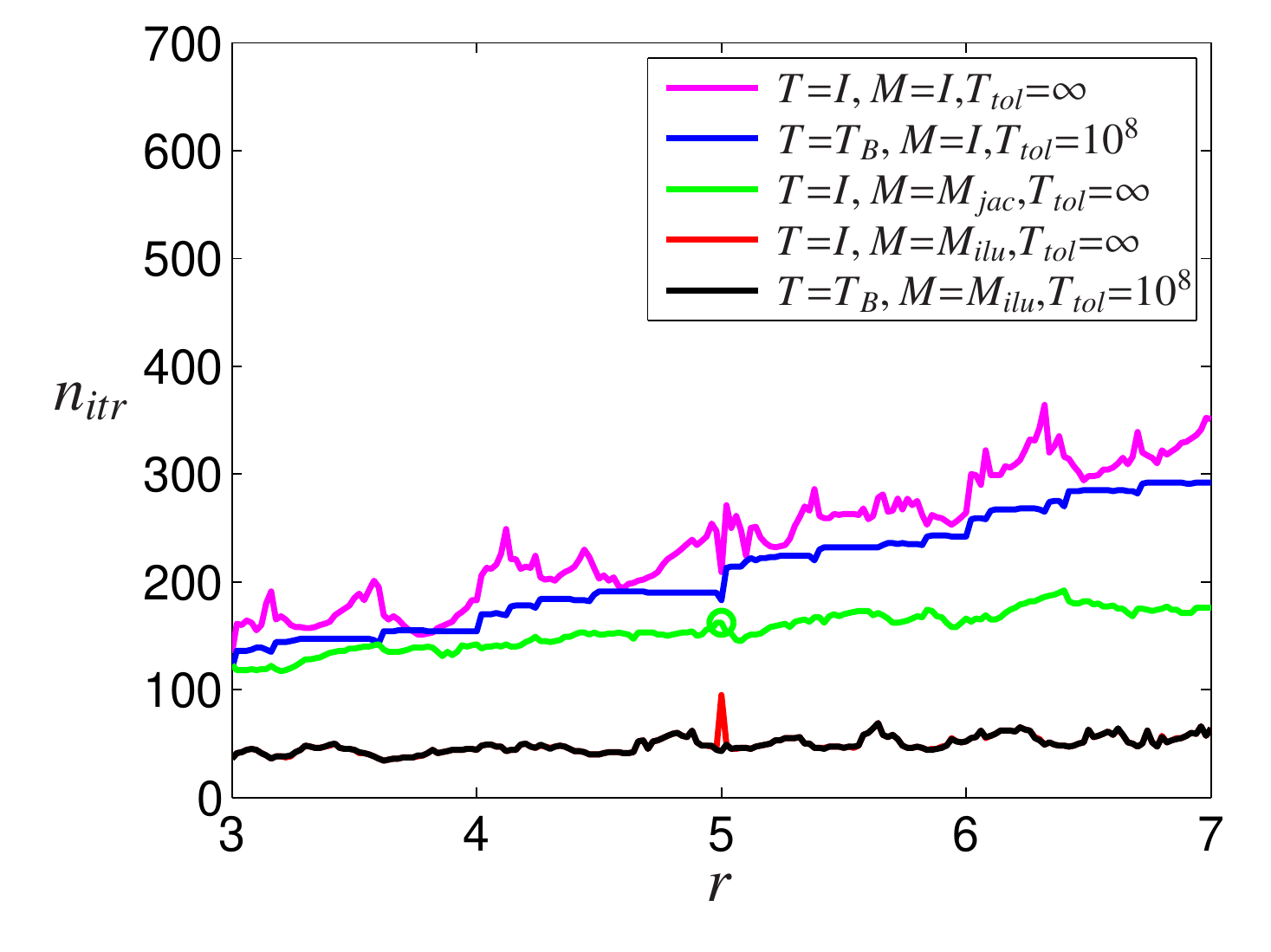}}
\subfloat[]{\includegraphics[trim=1.4in 5.2in 1in 2.1in,clip,width=0.33\textwidth]{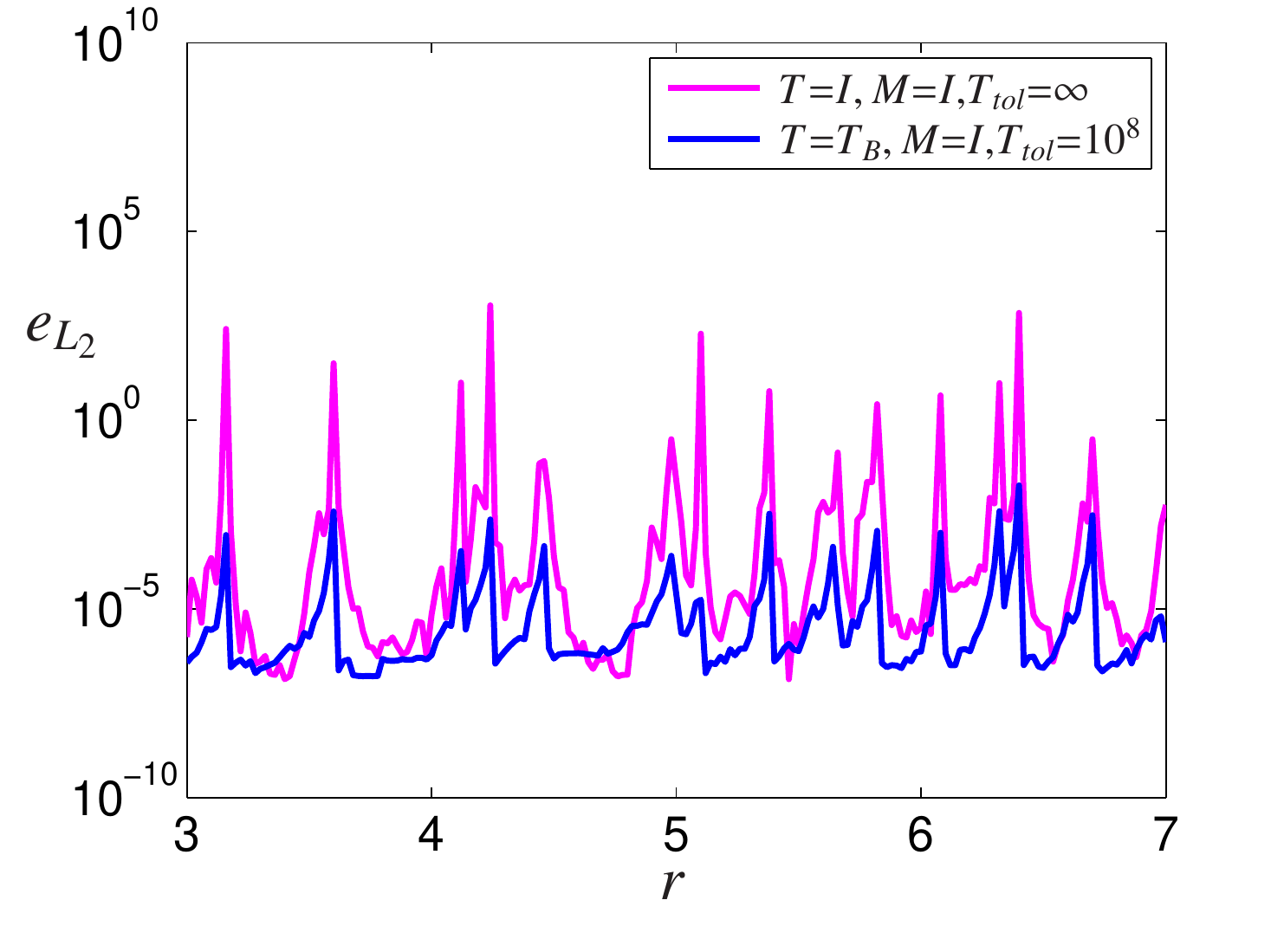}}
\subfloat[]{\includegraphics[trim=1.4in 5.2in 1in 2.1in,clip,width=0.33\textwidth]{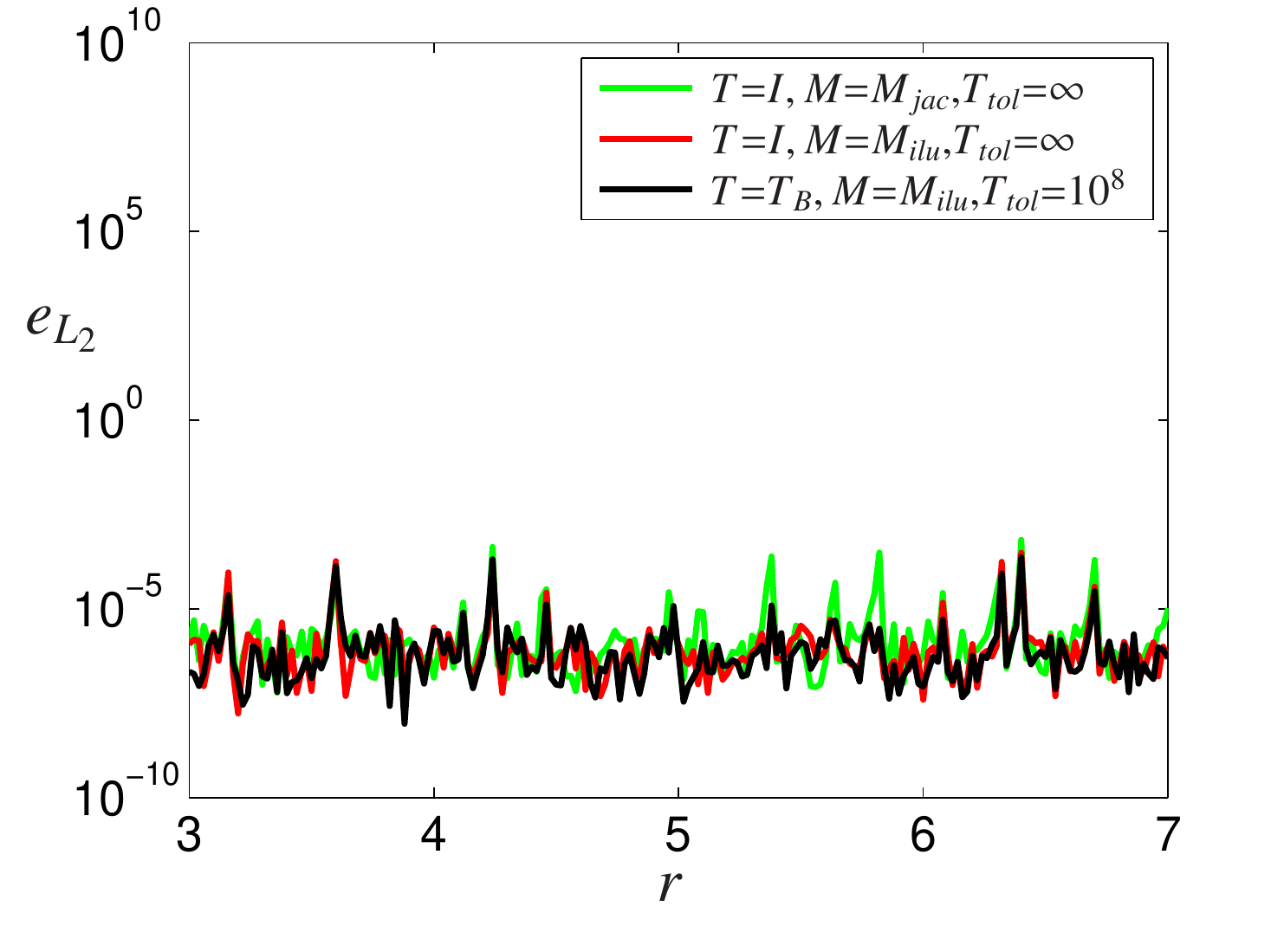}}
\caption{The (a) number of iterations, (b) $e_{L_2}$ with no solver preconditioner, and (c) $e_{L_2}$ with solver preconditioners $\bm{M}_{jac}$ and $\bm{M}_{ilu}$ using the stabilized Lagrange method. The open circles in (a) mark the values of $r$ at which the iterative solver failed to satisfy the stopping criteria.}
\label{fig:itrslag}
\end{figure*}

\begin{figure*}[ht]
\centering
\subfloat[]{\includegraphics[trim=1.4in 5.2in 1in 2.1in,clip,width=0.33\textwidth]{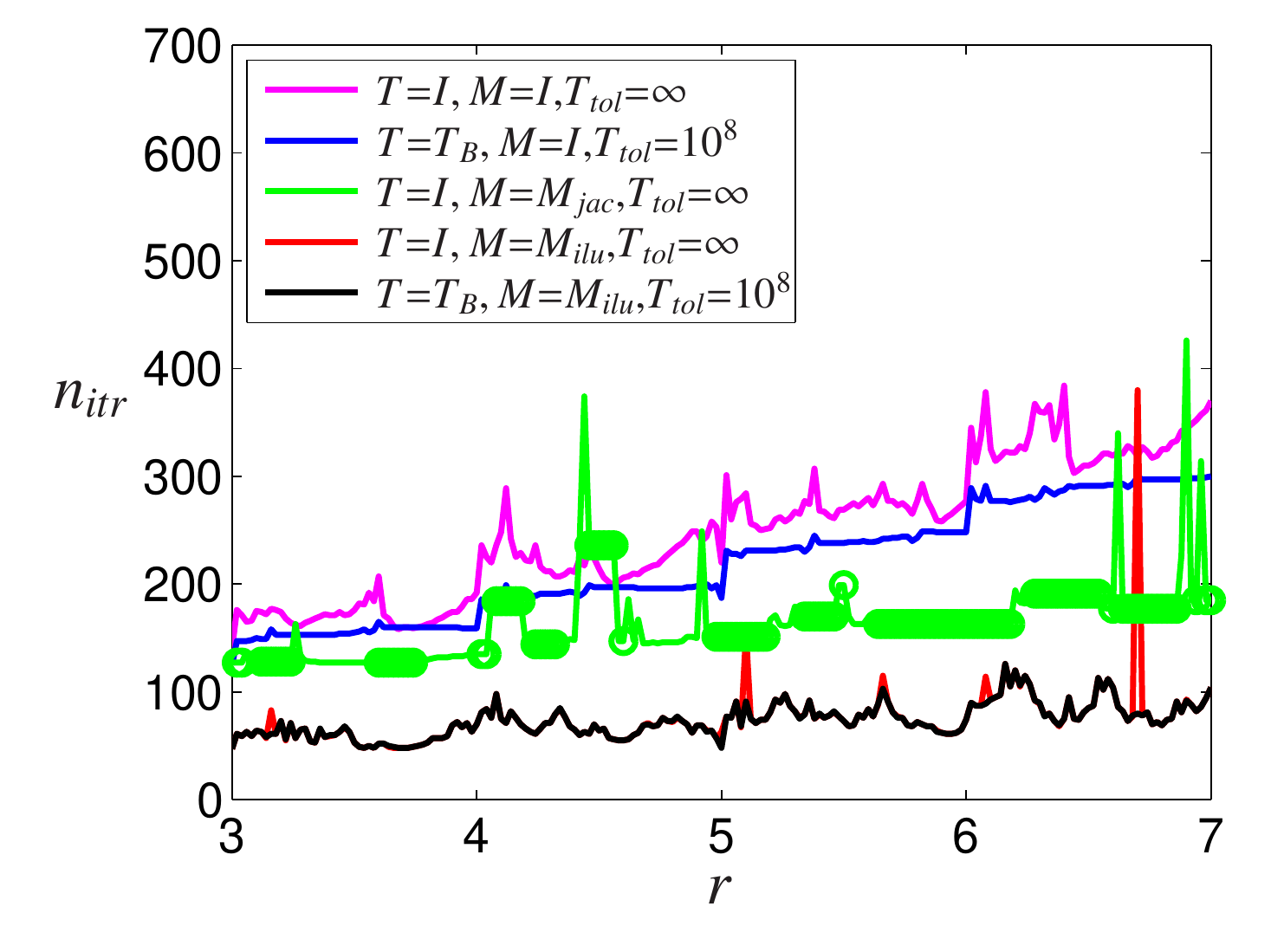}}
\subfloat[]{\includegraphics[trim=1.4in 5.2in 1in 2.1in,clip,width=0.33\textwidth]{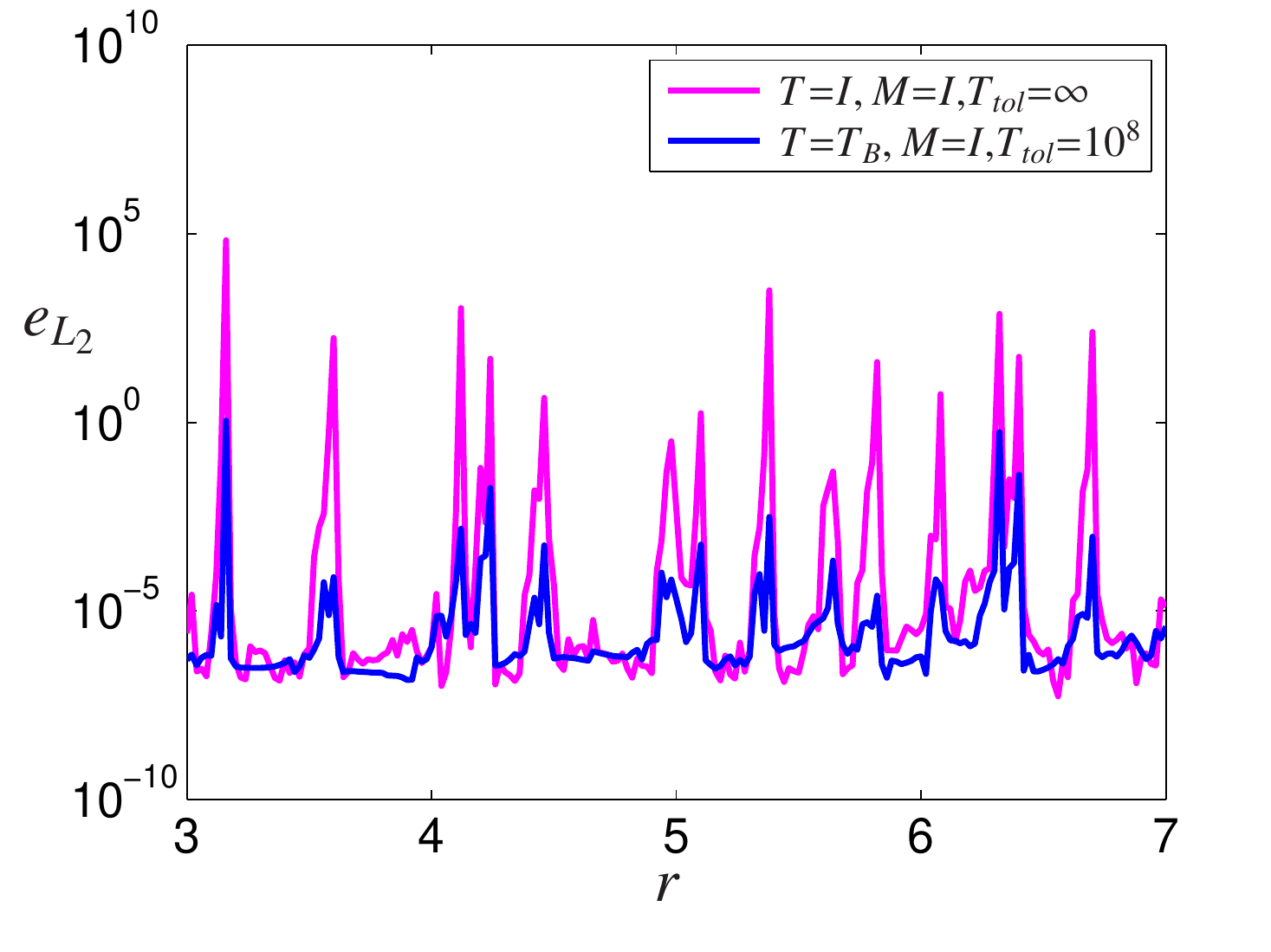}}
\subfloat[]{\includegraphics[trim=1.4in 5.2in 1in 2.1in,clip,width=0.33\textwidth]{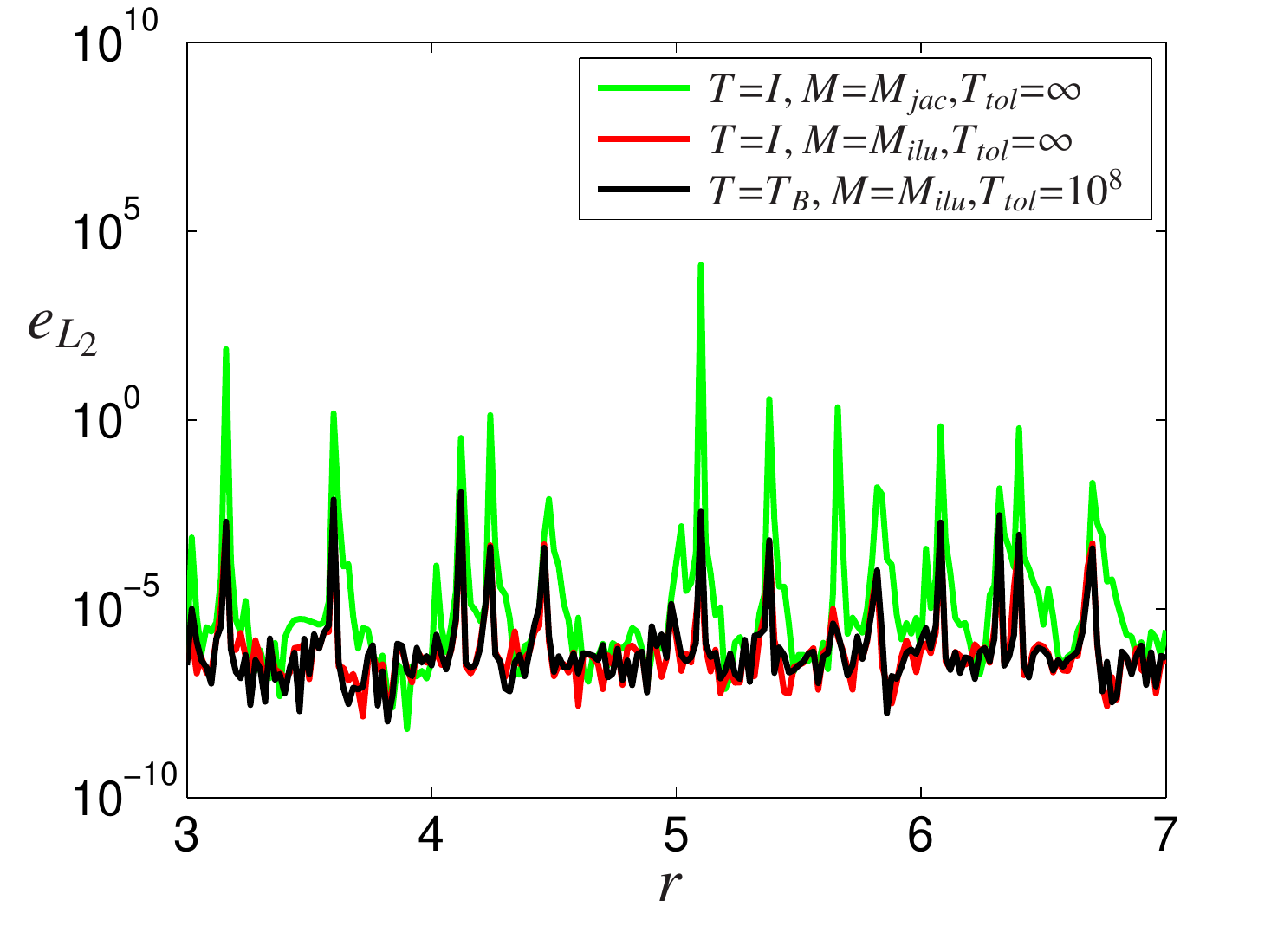}}
\caption{The (a) number of iterations, (b) $e_{L_2}$ with no solver preconditioner, and (c) $e_{L_2}$ with solver preconditioners $\bm{M}_{jac}$ and $\bm{M}_{ilu}$ using Nitsche's method. The open circles in (a) mark the values of $r$ at which the iterative solver failed to satisfy the stopping criteria.}
\label{fig:itrnits}
\end{figure*}

The third study examines the influence of the preconditioning scheme on the performance of an iterative solver by solving the system of equations using the generalized minimal residual method (GMRES) \citep{ref:saad1986}. A solver preconditioner, $\bm{M}$, was implemented to compare with the performance of the preconditioning scheme. A Jacobi, $\bm{M}_{jac}$, and incomplete LU with zero fill-in, $\bm{M}_{ilu}$, were chosen as the solver preconditioners. The number of iterations, $n_{itr}$, required to satisfy $\| \bm{f} - \bm{K}\hat{\bm{u}} \|_2 < 10^{-6}$ was determined using the physical solution. The solution error was determined as
\beq
e_{L_2} = \| \hat{u}-u_{ref2} \|_2 \text{ ,}
\eeq

\noindent where $u_{ref2}$ was a reference solution computed using a direct solver with $\bm{T}=\bm{I}$. The body-fitted FEM reference solution was not used here in order to distinguish the iterative solver error and the discretization error. Also, the reference solution at $r=5$ is not available because the direct solver fails due to the high condition number. Therefore $e_{L_2}$ is not computed at $r=5$. A comparison of the number of required iterations and the solution error is shown in Figs. \ref{fig:itrslag} and \ref{fig:itrnits} with and without $\bm{T}_B$ and $\bm{M}$. No degrees of freedom were constrained for $\bm{T}=\bm{I}$, denoted by $T_{tol}=\infty$.

As expected, the number of required iterations is reduced with the preconditioning scheme using both the stabilized Lagrange and Nitsche methods. The solver preconditioners reduce the number of iterations more than the projection scheme alone. However, the Jacobi preconditioner is not robust as the solver fails for some of the values of $r$. The geometric preconditioning scheme may be combined with a solver preconditioner. The incomplete LU preconditioner with and without $\bm{T}_B$ has the fewest required iterations. In this case, the proposed preconditioning scheme adds robustness, ensuring an almost constant number of iterations for all interface geometries.

\subsection{Example 3: Moving Cylinder in Channel Flow}

\begin{figure}[htb]
\centering\includegraphics[trim=2in 2.5in 2in 2.5in,clip,width=0.65\textwidth]{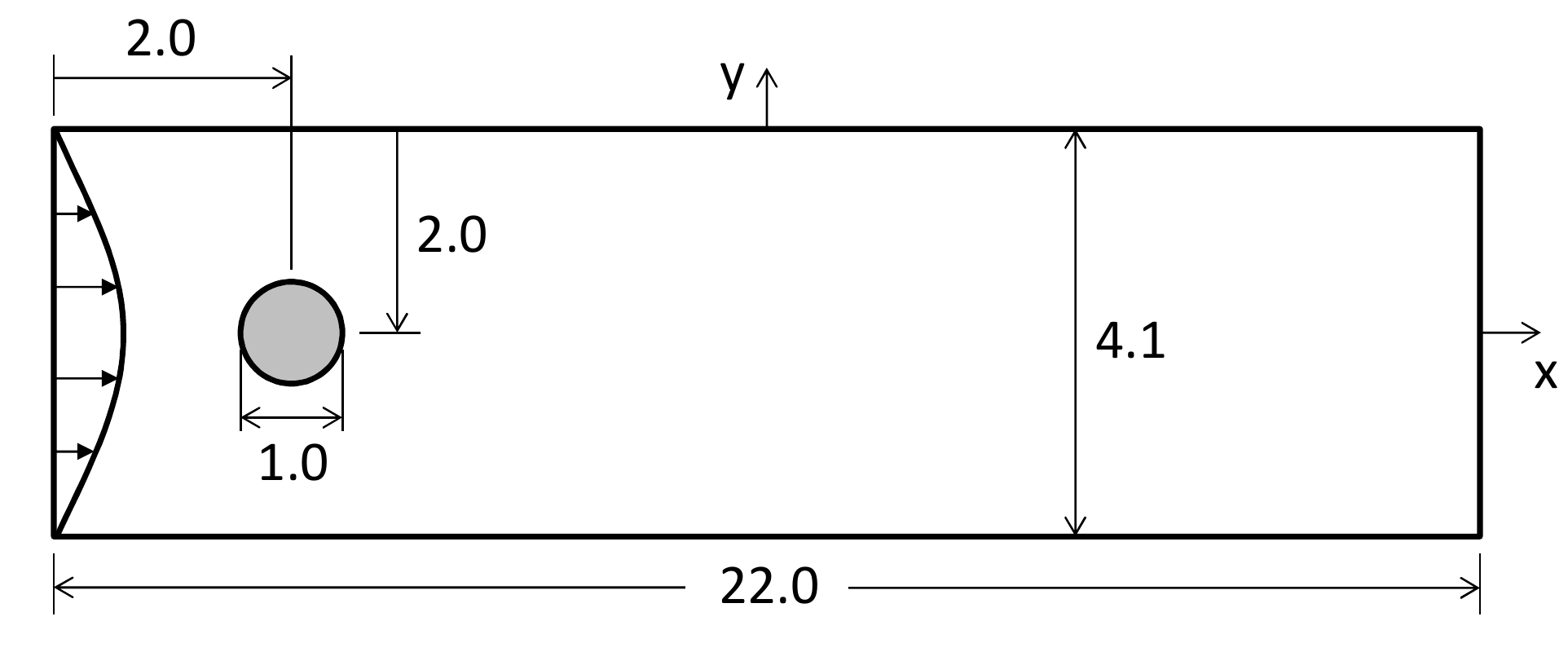}
\caption{Problem description for Example 3.}
\label{fig:channelgeometry}
\end{figure}

In this example, a 2D transient nonlinear problem with a moving interface is considered. A rigid cylinder immersed in a channel flow is oscillating perpendicular to the inflow direction. The flow is modeled by the incompressible Navier-Stokes equations, and the motion of the cylinder is prescribed by defining the level set field as an explicit function of time. The problem setup is depicted in Fig. \ref{fig:channelgeometry}. Note the fluid problem is modeled and solved in non-dimensional form. We study the stability and accuracy of the flow solution with and without the proposed preconditioning scheme for different $T_{tol}$ values for constraining degrees of freedom.

\begin{figure}[htb]
\centering\includegraphics[trim=1in 3.2in 1.5in 3.5in,clip,width=0.65\textwidth]{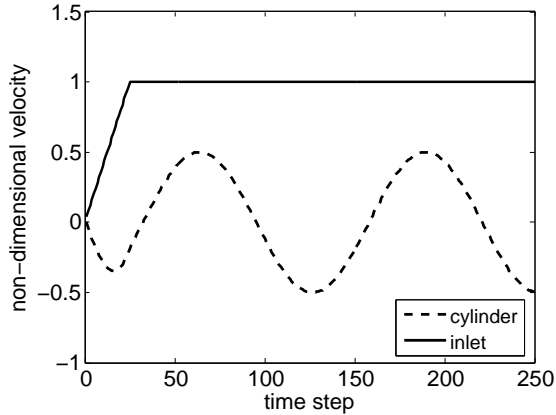}
\caption{Evolution of non-dimensional inlet and cylinder velocities.}
\label{fig:rampup}
\end{figure}

Along the channel inlet a parabolic inflow is prescribed. The outlet is assumed traction-free, and stick conditions are enforced at the upper and lower channel walls. The position of the cylinder and velocity along the cylinder surface, i.e. fluid-solid interface, are determined from the prescribed evolution of the discretized level set field.  The flow response is simulated over $250$ time steps with a non-dimensional time step size of $\Delta t = 0.05$. To facilitate the transient simulation of the flow field, we ramp up over time both the inlet conditions and the motion of the cylinder. The velocity profiles of the cylinder and inlet flow are depicted in Fig. \ref{fig:rampup}. The Reynolds number with respect to the maximum average inlet velocity is $100$. 

The weak form of the incompressible Navier-Stokes equations is discretized by four-node finite elements, i.e. the velocity and pressure fields are approximated piecewise by bilinear, equal-order interpolations. To avoid numerical instabilities we employ an SUPG/PSPG stabilization scheme \citep{TMRS:92}. The velocity boundary condition along the fluid-solid interface is enforced by a stabilized Lagrange multiplier formulation \citep{GeWa:10}. The Lagrange multipliers are approximated element-wise by bilinear shape functions. The reader is referred to Kreissl and Maute \citep{KM:12} for additional details on the XFEM implementation of the flow model. The flow solution is advanced in time with an Euler-backward time integration scheme. In each time step, the nonlinear sub-problem is solved by the Newton-Raphson method, and a direct solver is applied to the linearized problem. The nonlinear residual is required to drop by only $10\%$ in each time step.

\begin{figure}[htb]
\centering
\subfloat[]{\includegraphics[trim=1.25in 3.2in 1.6in 3.5in,clip,width=0.65\textwidth]{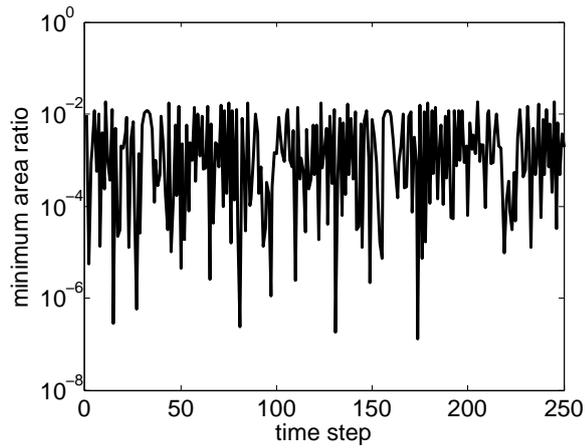}}\\
\subfloat[]{\includegraphics[trim=1.50in 5.2in 1.25in 2in,clip,width=0.65\textwidth]{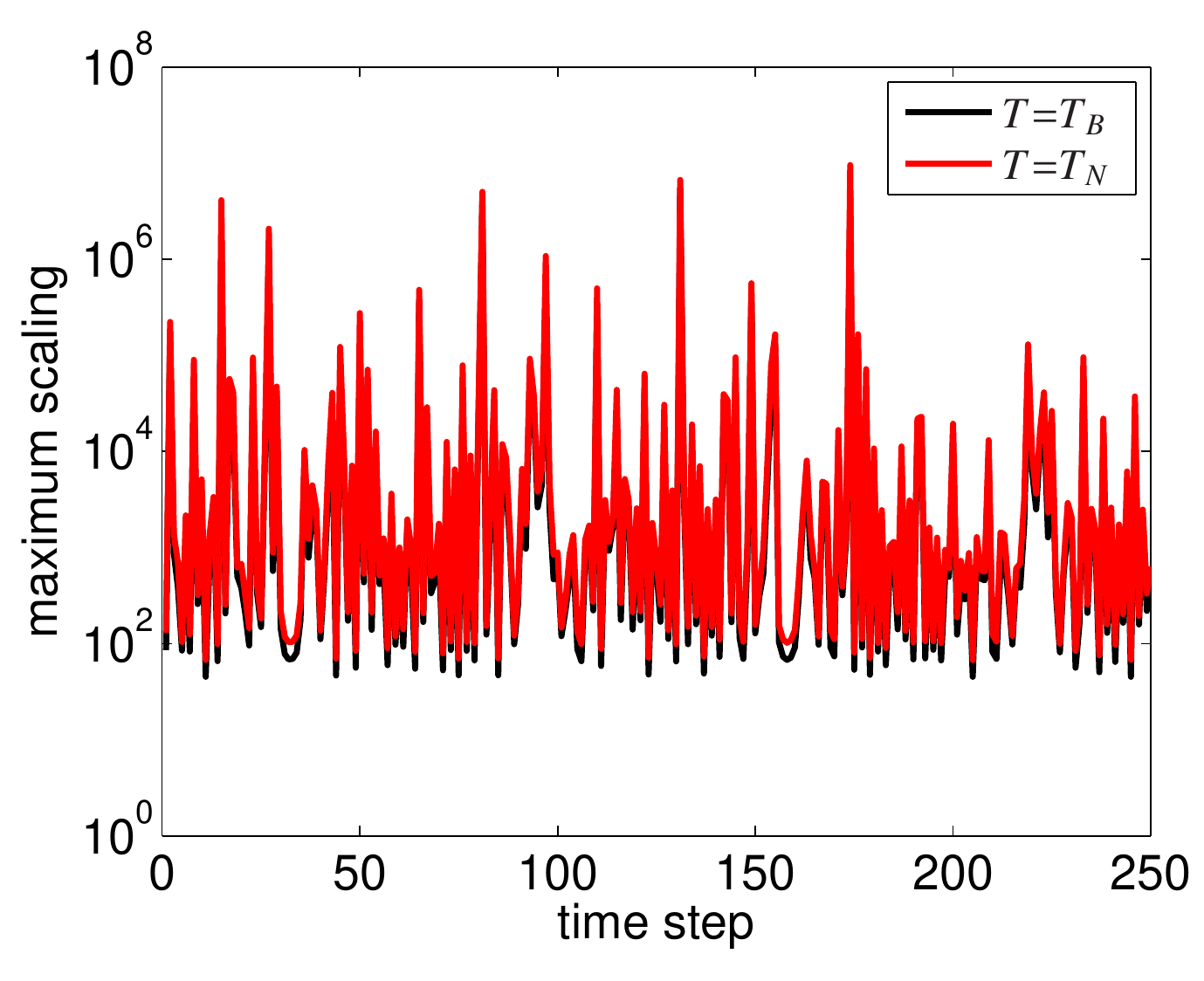}}
\caption{The (a) minimum elemental fluid area ratio, and (b) the maximum entry in the preconditioning matrices $\bm{T}_N$ and $\bm{T}_B$ in each time step.}
\label{fig:evolratios}
\end{figure}

First we discretize the channel with $6912$ elements and $7105$ nodes. The mesh in the vicinity of the cylinder is uniform with a non-dimensional element size of $0.085 \times 0.085$. Initially, $48$ elements are intersected by the fluid-solid interface and the flow field is approximated by $63,255$ degrees of freedom. As the cylinder oscillates, the intersection configuration, number of intersected elements, and number of degrees of freedom change slightly. The evolutions of the minimum ratio of elemental fluid area over the total elemental area and the maximum entries in the preconditioning matrix for the $\bm{T}_N$ and $\bm{T}_B$ formulations are shown in Fig. \ref{fig:evolratios}. The minimum area ratios of order $10^{-7}$ lead to large entries in the preconditioning matrix. The maximum entries in $\bm{T}_N$ are slightly larger than the ones in $\bm{T}_B$ but are of the same order. The evolution of both formulations is similar.

We compare the performance of the proposed preconditioning scheme against an approach where only degrees of freedom with vanishing influence are constrained. We examine the evolution of the total horizontal and vertical forces acting on the cylinder, and we consider $T_{tol}=[10^8,10^6,10^4,10^2,10^1]$ for determining the constrained degrees of freedom. As shown in Fig. \ref{fig:evolratios}, the maximum value of the preconditioning matrix is less than $10^8$ for all time steps. Therefore no degrees of freedom are constrained for $T_{tol}=10^{8}$. For $T_{tol}\le10^6$, the number of constrained degrees of freedom increases as $T_{tol}$ is reduced. The number of constrained degrees of freedom varies with time, and the maximum is shown in Table \ref{tab:dropdofs} when $\bm{T}_B$ is applied.

\begingroup
\renewcommand*{\arraystretch}{1.2}
\begin{table}[ht]
  \centering
    \begin{tabular}{cc}
    \hline
    $T_{tol}$ & max constrained dofs \\
    \hline
    $10^8$ & $0$  \\
    $10^6$ & $3$  \\
    $10^4$ & $9$  \\
    $10^2$ & $18$ \\
    $10^1$ & $30$ \\
    \hline
    \end{tabular}
  \caption{Maximum number of constrained degrees of freedom.}
  \label{tab:dropdofs}
\end{table}
\endgroup

\begin{figure}[htb]
\centering
\subfloat[]{\includegraphics[trim=1.25in 5.2in 1.25in 2in,clip,width=0.65\textwidth]{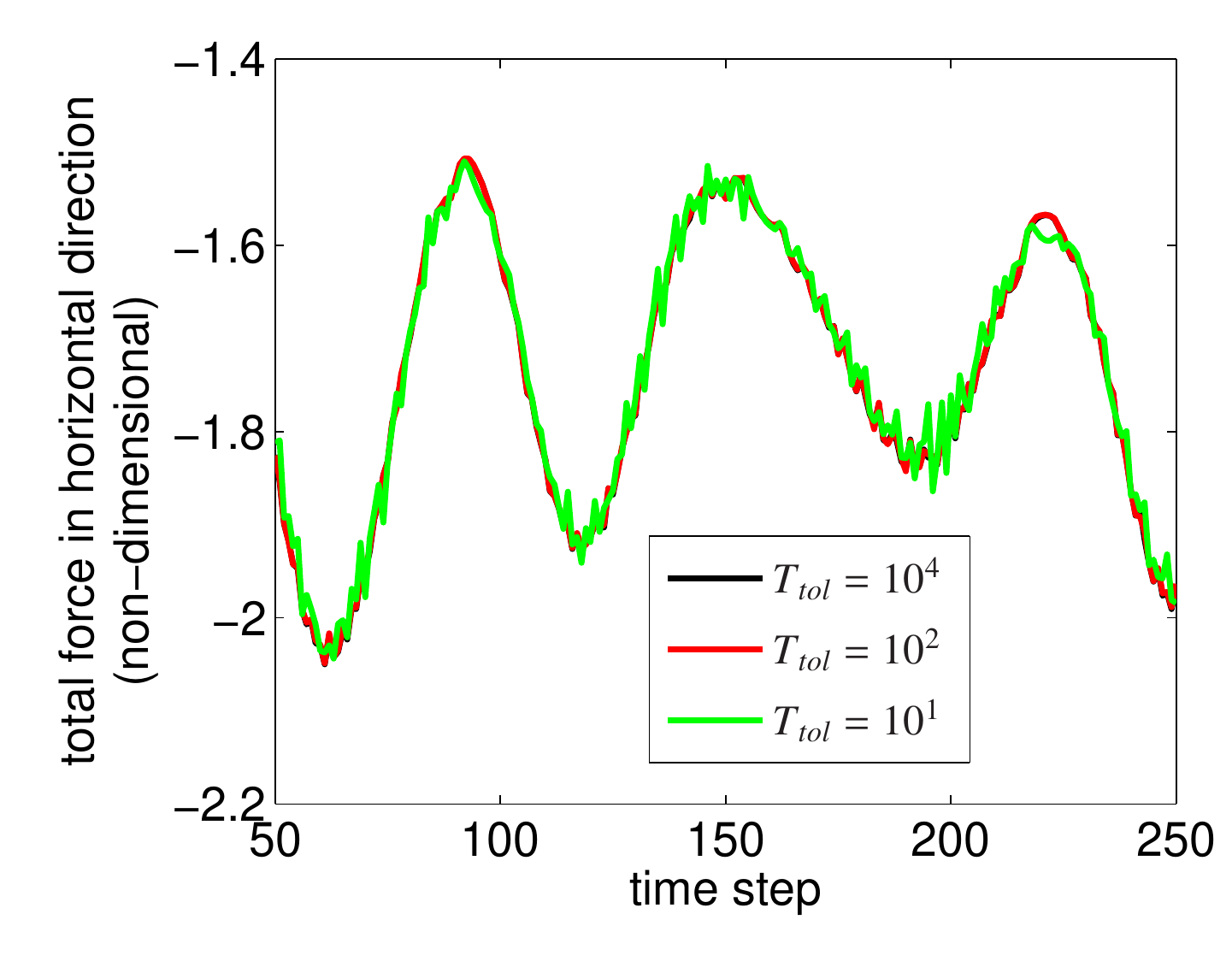}}\\
\subfloat[]{\includegraphics[trim=1.25in 5.2in 1.25in 2in,clip,width=0.65\textwidth]{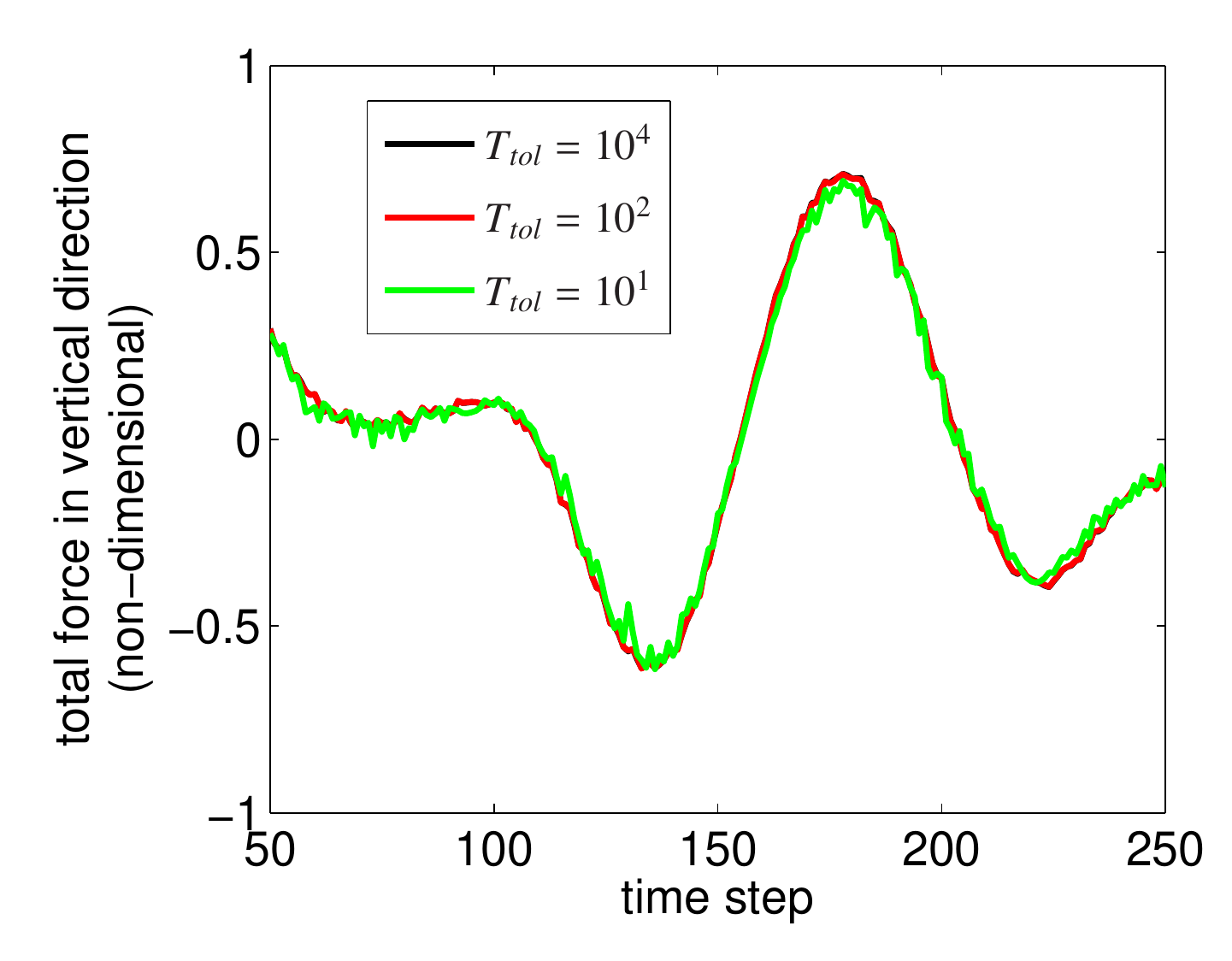}}
\caption{Evolution of the total force in the (a) horizontal and (b) vertical directions for different values of $T_{tol}$ when no geometric preconditioner is applied.}

\label{fig:noprec}
\end{figure}

Without the proposed preconditioning scheme, the transient simulation diverges for $T_{tol}>10^{4}$. The evolutions of the total horizontal and vertical forces for $T_{tol}=[10^4,10^2,10^1]$ are depicted in Fig. \ref{fig:noprec}. Note, the results are shown only after $50$ time steps for which the influence of ramping up the inlet and cylinder velocities has sufficiently faded. The force evolutions for $T_{tol}=10^4$ and $T_{tol}=10^2$ are similar. However, if $T_{tol}$ is chosen too low, here $T_{tol}=10^1$, the forces erroneously oscillate. The proper choice of $T_{tol}$ is not known {\it a priori}. 

\begin{figure}[htb]
\centering
\subfloat[]{\includegraphics[trim=1.25in 5.2in 1.25in 2in,clip,width=0.65\textwidth]{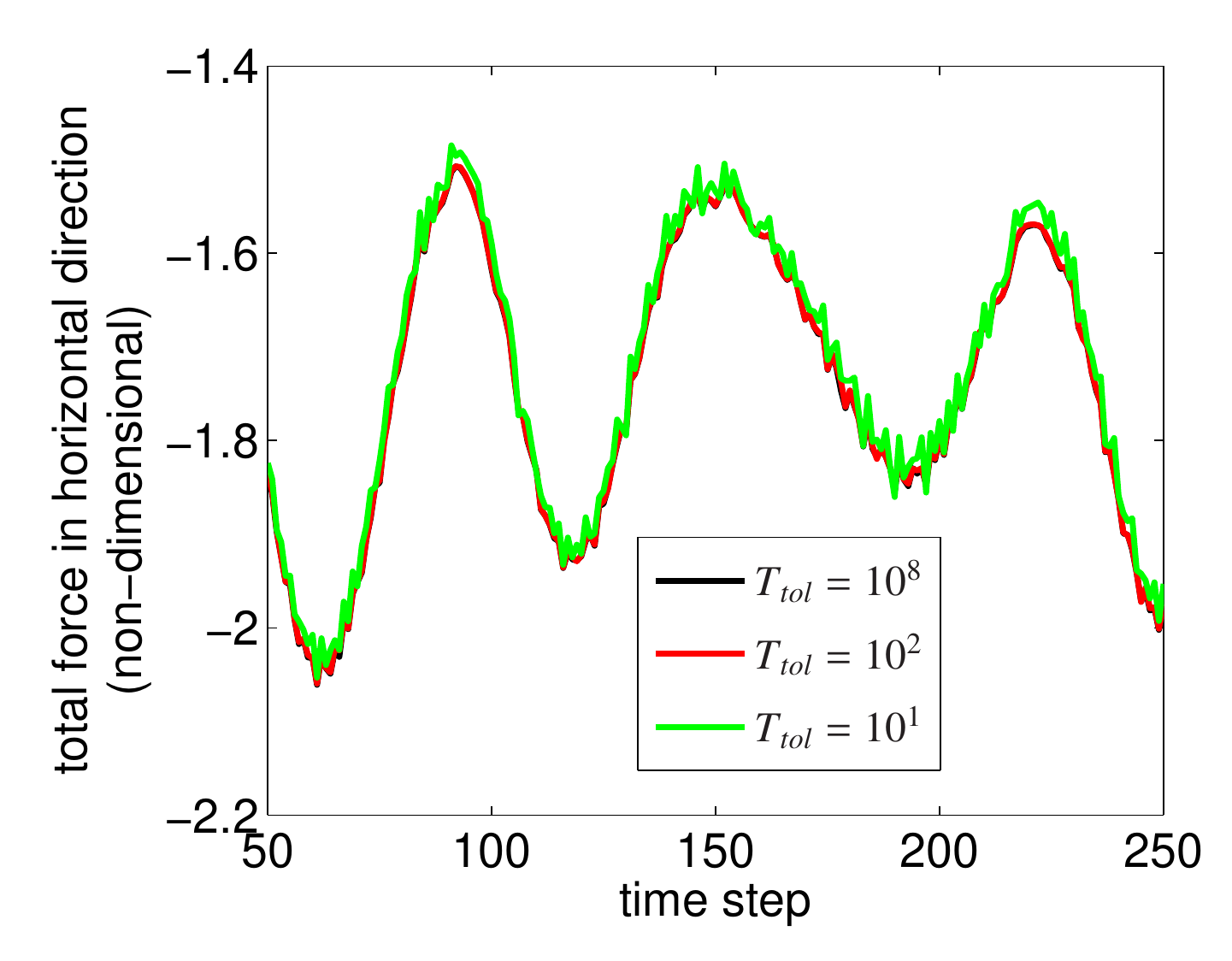}}\\
\subfloat[]{\includegraphics[trim=1.25in 5.2in 1.25in 2in,clip,width=0.65\textwidth]{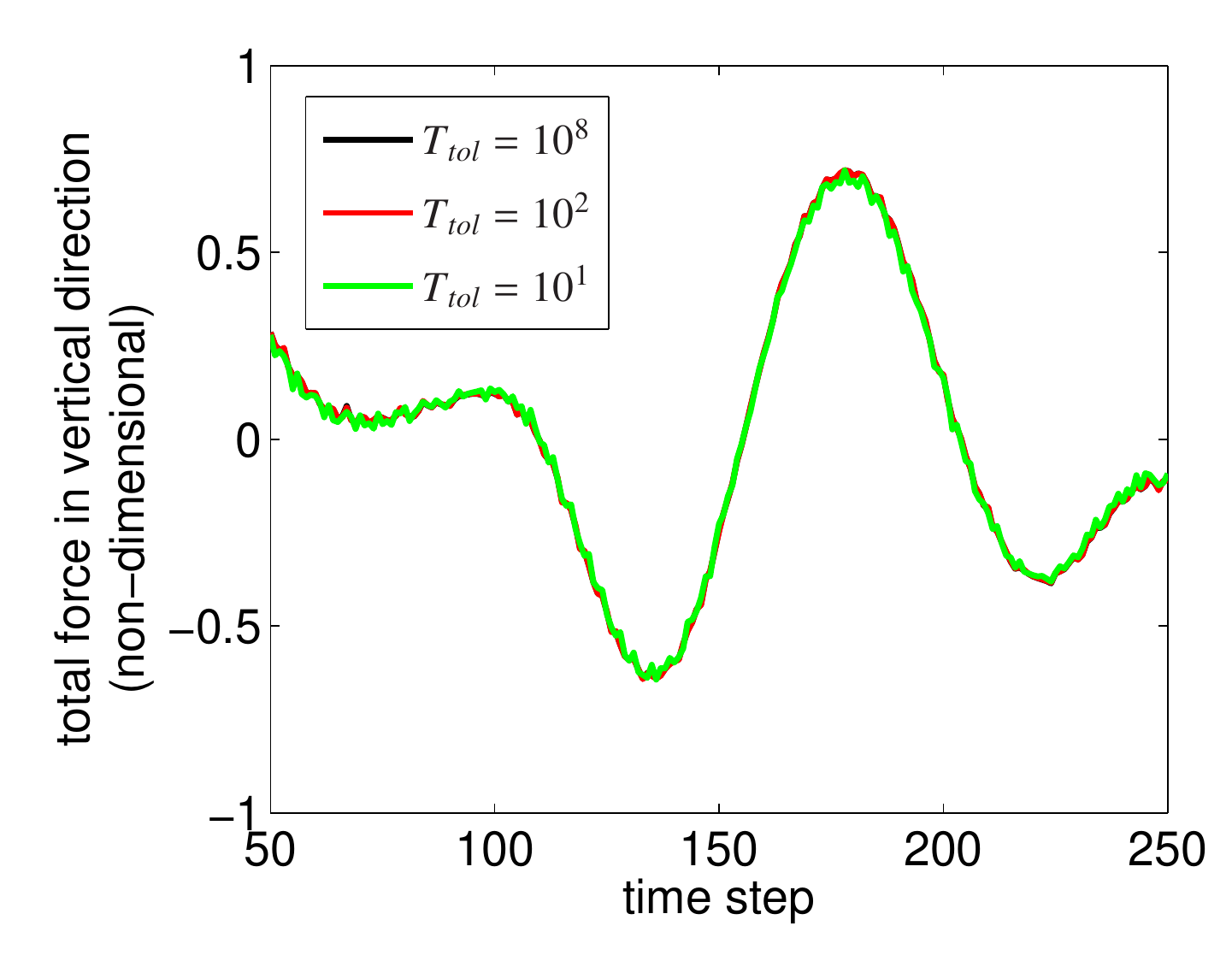}}
\caption{Evolution of the total force in the (a) horizontal and (b) vertical directions for different values of $T_{tol}$ using the $\bm{T}_B$ preconditioning matrix.}
\label{fig:bprec}
\end{figure}

In contrast, no convergence issues were observed with the proposed preconditioning scheme for both formulations of the preconditioning matrix. In Fig. \ref{fig:bprec}, the evolution of the total forces are shown using the $\bm{T}_B$ preconditioning matrix. Note, as the results for $T_{tol}=[10^8,10^6,10^4]$ are indistinguishable, only the results for $T_{tol}=[10^8,10^2,10^1]$ are shown. The results for the $\bm{T}_N$ preconditioning matrix are nearly identical and therefore not shown. For both formulations, a similar behavior can be observed when no preconditioner is used: if $T_{tol}$ is too low the forces erroneously oscillate with a high frequency. This behavior seems to be less pronounced when the preconditioning matrix is used. 

\begin{figure}[htb]
\centering\includegraphics[trim=1.25in 5.2in 1.25in 2in,clip,width=0.65\textwidth]{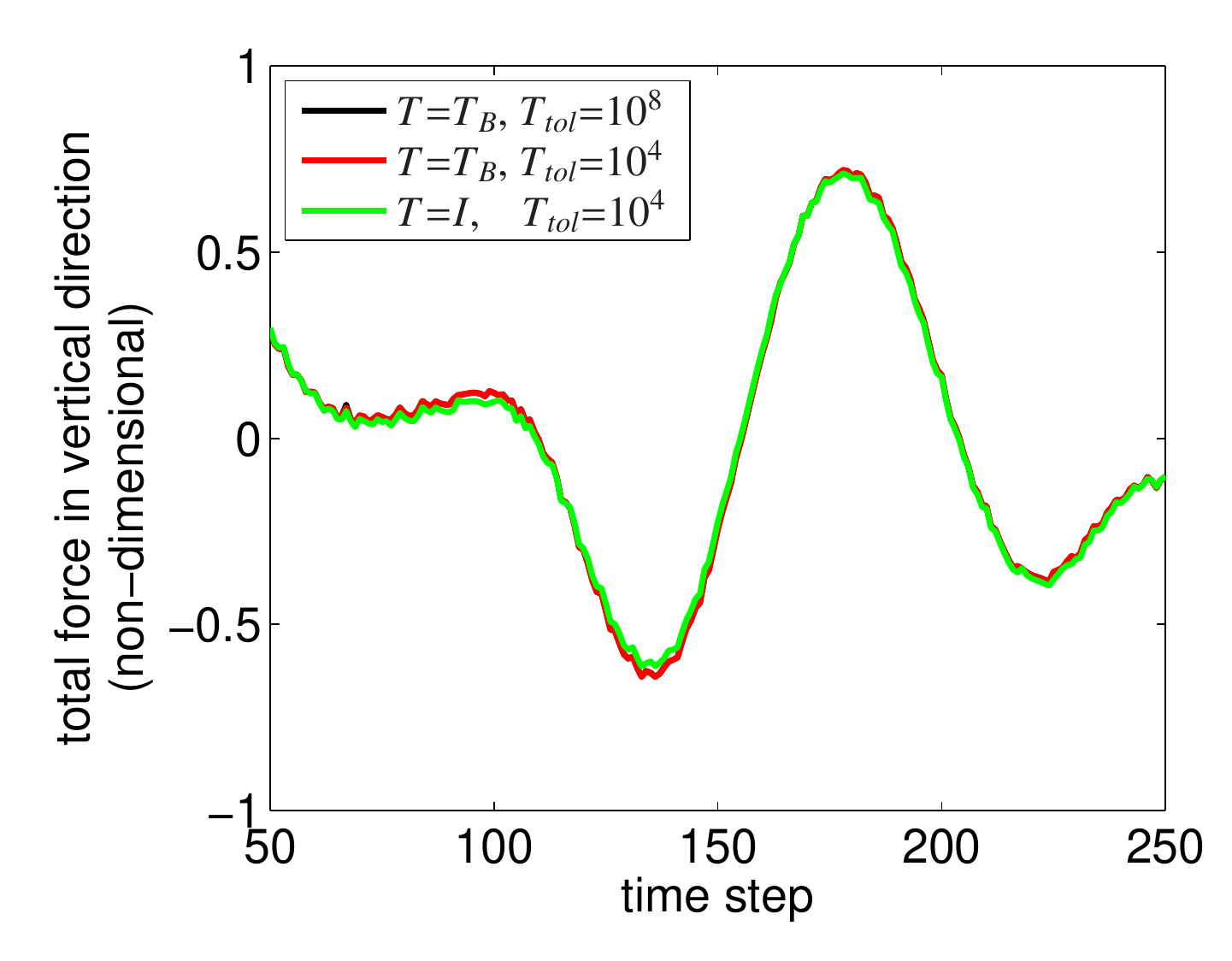}
\caption{Comparison of the total vertical force for different values of $T_{tol}$ with and without $\bm{T}_B$.}
\label{fig:comp_frcy}
\end{figure}

A direct comparison of the results obtained with and without the $\bm{T}_B$ preconditioning matrix are depicted in Fig. \ref{fig:comp_frcy}. Here we only compare the total force in vertical direction for $T_{tol}=[10^8,10^4]$. The results for the total force in the horizontal direction show similar behaviors and are therefore omitted. Recall the simulations diverge for $T_{tol}=10^8$ when no preconditioner is used. While the results with $\bm{T}_B$ are indistinguishable for $T_{tol}=[10^8,10^4]$, the cross-comparison between the force evolutions for $T_{tol}=10^4$ shows a slight difference. This is attributed to the different convergence behavior; the convergence of the Newton-Raphson method is once monitored in the physical and once in the transformed space. As a stricter convergence is enforced, the difference decreases. 

\begin{figure}[htb]
\centering
\subfloat[]{\includegraphics[trim=1in 3.2in 1.5in 3.5in,clip,width=0.65\textwidth]{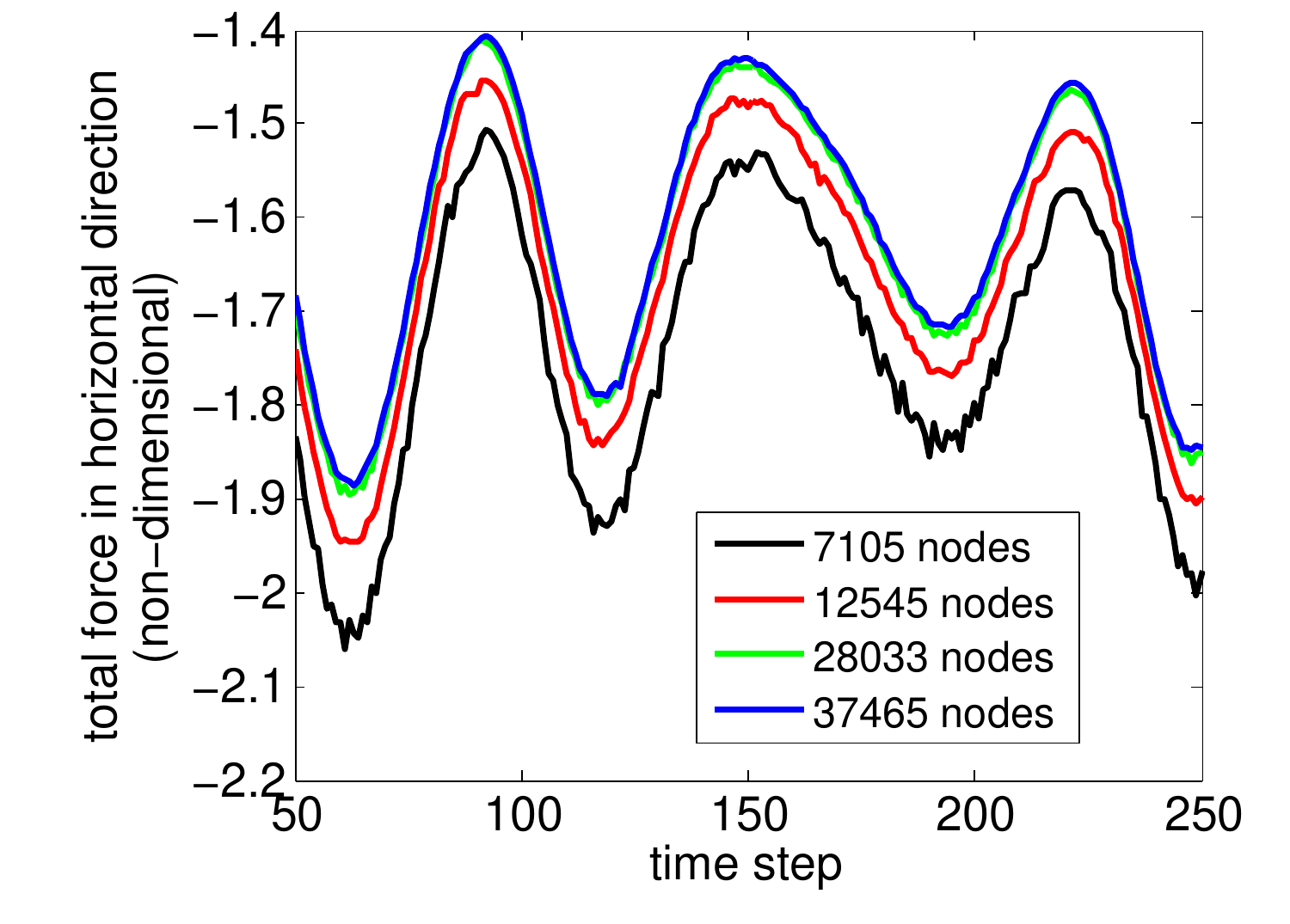}}\\
\subfloat[]{\includegraphics[trim=1in 3.2in 1.5in 3.5in,clip,width=0.65\textwidth]{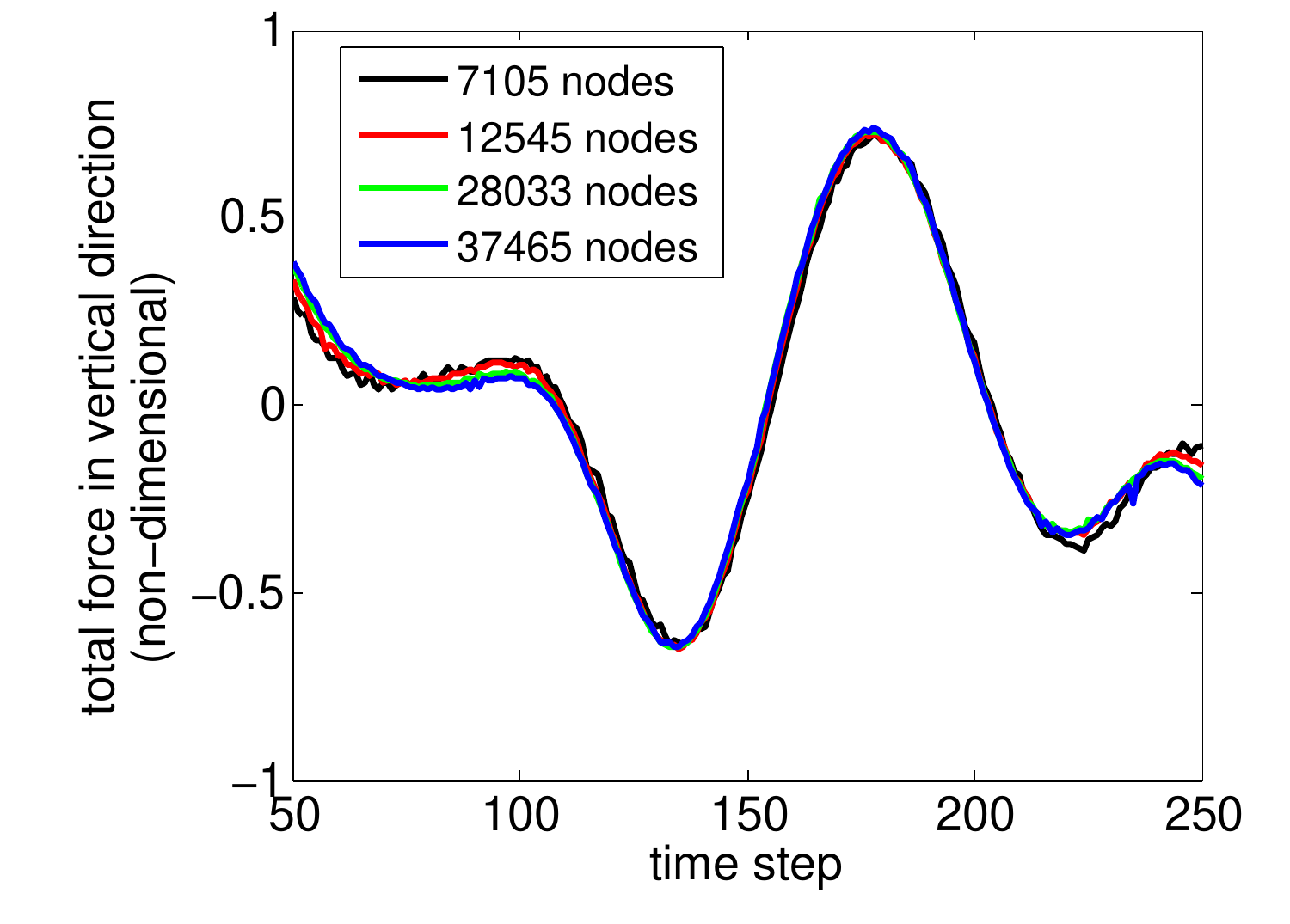}}
\caption{Evolution of total vertical forces in horizontal and vertical directions for different mesh refinement levels using the $\bm{T}_B$ formulation of the geometric preconditioner with $T_{tol}=\infty$.}
\label{fig:meshref}
\end{figure}

The robustness provided by the preconditioning \linebreak scheme allows the problem to be solved on refined \linebreak meshes. We examine the total horizontal and vertical forces acting on the cylinder using different mesh sizes. The considered mesh sizes and the number of initially intersected elements are given in Table \ref{tab:meshref}. The evolution of the total force for a sequence of refined meshes is shown in Fig. \ref{fig:meshref}. For $T_{tol}=10^{8}$, no degrees of freedom were constrained for all mesh sizes. The force evolutions converge as the mesh is refined. The high frequency oscillations present in the coarsest mesh vanish with mesh refinement.

\begin{table}[htb]
	\centering
		\begin{tabular}{ccc}
		\hline
		nodes & elements & intersected elements \\
		\hline
		$7105$  & $6912$  & $48$ \\
		$12545$ & $12288$ & $64$ \\
		$28033$ & $27648$ & $92$ \\
		$37465$ & $37044$ & $120$  \\
		\hline
		\end{tabular}
	\caption{Mesh refinement study.}
	\label{tab:meshref}
\end{table}

\section{Conclusions}

A simple and efficient preconditioning scheme has been proposed for Heaviside enriched XFEM problems which transforms the discretized governing equations into a well-conditioned form. The preconditioning scheme consists of a geometric preconditioner and constraining degrees of freedom to zero which interpolate the solution for small areas of intersection. The geometric preconditioner is constructed from the nodal basis functions and the interface configuration. Therefore the preconditioning matrix can be computed prior to constructing the system matrices, making it well-suited for nonlinear problems. The ill-conditioning due to small element intersections is eliminated, and the condition number of the system matrices is comparable to that of a body-fitted mesh using the traditional FEM. 

We have shown that when only selecting degrees of freedom to constrain to zero without the preconditioning matrix, there is a strong trade-off with reducing the condition number and a loss in solution accuracy. By implementing the proposed preconditioning scheme the condition number is reduced, and a loss in solution accuracy only occurs if the tolerance criteria for selecting the degrees of freedom to constrain is too small. While generic solver preconditioners help to reduce the condition number, the proposed preconditioning scheme is robust and efficient for solving linear and nonlinear problems. Additionally, the proposed approach performs well for the stabilized Lagrange and Nitsche methods for enforcing continuity at the interface.

In this work two diagonal forms of the preconditioning matrix were studied. Additional approaches for building the preconditioning matrix can be further explored, including diagonal and non-diagonal forms. Only 2D problems with static and prescribed moving interfaces were considered. The extension of the proposed preconditioning scheme to 3D problems is straight forward. The performance of the preconditioning scheme for problems with dynamically evolving interfaces will be investigated in future studies. 

\section*{Acknowledgments}
The first author acknowledges the support of the NASA Fundamental Aeronautics Program Fixed Wing Project, and the second and fourth authors acknowledges the support of the National Science Foundation under grant CMMI-0729520. The third author acknowledges the support of the Department of Energy under grant DE-SC0006402. The opinions and conclusions presented are those of the authors and do not necessarily reflect the views of the sponsoring organizations.

\bibliographystyle{elsarticle-num}
\bibliography{refs}		

\end{document}